\documentclass[11pt]{article}

\usepackage{amssymb, epsfig, amsmath,amsfonts, subfigure}
\usepackage{hyperref}

\usepackage{rotating}
\usepackage[toc,page]{appendix}
\usepackage{color}
\usepackage{multirow}

\usepackage{subfigure}
\usepackage{graphics}
\usepackage{graphicx}
\usepackage{color}
\usepackage{soul}
\usepackage[makeroom]{cancel}

\usepackage{latexsym}
\usepackage{graphics}
\usepackage{amsmath}
\usepackage{amsthm}
\usepackage{xspace}
\usepackage{amssymb}
\usepackage{color}
\usepackage{cite}
\usepackage{latexsym}
\usepackage{graphics}
\usepackage{amsmath}
\usepackage{amsthm}
\usepackage{xspace}
\usepackage{amssymb}
\usepackage{epsfig}
\usepackage{comment}

\newcommand{\beql}[1]{\begin{equation}\label{#1}}
\newcommand{\eeql}{\end{equation}}
\newcommand{\eqn}[1]{(\ref{#1})}

\newcommand{\R}{\mathbb{R}}
\newcommand{\pr}{\mathbb{P}}
\newcommand{\E}{\mathbb{E}}

\newtheorem{thm}{Theorem}
\newtheorem{lem}[thm]{Lemma}

\newtheorem{cor}[thm]{Corollary}

\newcommand{\RR}{{\mathbb R}}

\newcommand{\NN}{{\mathbb N}}

\newcommand{\ra}{\rightarrow}

\newcommand{\eeq}{\end{equation}}
\newcommand{\beqal}[1]{\begin{eqnarray}\label{#1}}
\newcommand{\eeqa}{\end{eqnarray}}
\newcommand{\eeqno}{\end{displaymath}}

\newcommand{\la}{\lambda}
\newcommand{\ep}{\epsilon}

\newcommand{\bone}{{\mathbf 1}}

\newcommand{\qandq}{\quad\mbox{and}\quad}

\newcommand{\qasq}{\quad\mbox{as}\quad}

\newcommand{\non}{\nonumber}
\newcommand{\RA}{\Rightarrow}
\DeclareMathOperator{\sgn}{sgn}

\begin{document}

\title{A service system with on-demand agent invitations
}

\author
{Guodong Pang\\
Harold and Inge Marcus Department\\
of Industrial and Manufacturing Engineering\\
Pennsylvania State University\\
University Park, PA 16802\\
\texttt{gup3@psu.edu}
\and
Alexander L. Stolyar \\
Department of Industrial\\ and Systems Engineering\\
Lehigh University\\
Bethlehem, PA 18015 \\
\texttt{stolyar@lehigh.edu}
}

\date{\today}

\maketitle

\allowdisplaybreaks

\begin{abstract}

We consider a service system where agents are invited on-demand.  Customers arrive exogenously as a Poisson process and join a customer queue upon arrival if no agent is available.   Agents decide to accept or decline invitations after some exponentially distributed random time, and join an agent queue upon invitation acceptance if no customer is waiting. 
A customer and an agent are matched in the order of customer arrival and agent invitation acceptance under the non-idling condition, and will leave the system simultaneously once matched (service times are irrelevant here). 

We consider a feedback-based adaptive agent invitation scheme, which controls 
the number of pending agent invitations,
depending on the customer and/or agent queue lengths and their changes.
The system process has two components -- `the difference between agent and customer queues' and `the number of pending invitations', and is
a countable continuous-time Markov chain.

For the case when the customer arrival rate is constant, we establish fluid and diffusion limits, in the asymptotic regime
where the customer arrival rate goes to infinity, while the agent response rate is fixed.
We prove the process stability and fluid-scale limit interchange, which in particular imply that both customer and
agent waiting times in steady-state vanish in the asymptotic limit. To do this we develop a novel (multi-scale)  Lyapunov drift argument;
it is required because the process has non-trivial behavior on the state space boundary.
When the customer arrival rate is time-varying, we present a fluid limit for the processes in the same asymptotic
regime. Simulation experiments are conducted to show
good performance of the invitation scheme and accuracy
of fluid limit approximations.

\end{abstract}


\section{Introduction}
\label{sec-intro}

The model in this paper is primarily motivated by applications to call/contact
 centers, where customer requests (arriving as calls, chat sessions, etc.),
are answered and processed by agents.  
(The model has broader potential applications,  for example, to inventory models 
as discussed later.) 
In such systems,
some of the callers/customers need to be served not by regular agents, 
who are always available (either working on a call or standing by)
during their work hours, but different
agents who are highly skilled and/or have knowledge that typical agents may not possess. These agents with special expertise/knowledge, are referred to as {\it special agents}, or {\em knowledge workers} \cite{Patent2010}. 
The time of knowledge workers is usually very valuable (and expensive).
As a result, it is not rational, or even feasible, to have a pool of
knowledge workers constantly available. Instead, they are invited on-demand
(usually via remote online access). This mode of operation, however, poses
a challenge, because (a) knowledge workers that accepted an invitation
should not wait for an extended (expensive) time before they actually start 
processing a call,
and (b) the callers should not have long waiting time either,
because some of the valuable calls will be lost and the service level objectives of the call center will not be achieved. Thus, effective mechanisms must be designed to assure that both these objectives are achieved.

Therefore, a critical part of the operation procedure for such a system,
is an algorithm 
to decide when to invite knowledge workers and how many to invite in order to minimize/stabilize both customer and agent delays. (From now on we will often say ``agent" instead of ``knowledge worker", as this is the only agent type that we consider
in this paper.) 
One ``naive" scheme would be to invite an agent for each arrived customer -- this is grossly inefficient, because 
an agent responds after (often large) delay, and may not respond at all. 
Another ``naive" approach
could be that agents are invited  at the rate that is equal to the rate of customers arrivals. However, {\em even if the agent arrival rate is known},
 it can be easily checked that the variance of the difference of the number of
customer arrivals and agents invitations increases linearly as time evolves.
Let alone the fact that the agent arrival rate is usually not known exactly 
in advance and may change over time.
Thus, this scheme is not desirable.

A feedback invitation scheme is proposed for the management of knowledge workers in \cite{Patent2010}. The idea is to keep a dynamic target level of invited agents, which changes according  
to the current system state and its ``derivative''.
The scheme allows to stabilize the system and keep
waiting times of customers and agents low.
In this paper we consider a stylized model of that invitation scheme -- it 
is such that the stochastic process describing
 the system dynamics is a continuous-time Markov chain (CTMC). 
Simulation experiments confirm that our stylized model's behavior is
 indeed very close to that of the (more practical)
feedback invitation scheme in \cite{Patent2010}; see section \ref{secNumerical}. 
Our analysis of the model shows that waiting times of both customers  and agents
are stable and
asymptotically negligible, as the system scale (the customer
arrival rate) goes to infinity.

Specifically, in the stylized model, 
agents are invited to serve customers, and  decide to accept or reject the invitations after some random delay. 
Customers join a customer queue upon arrival if no agent is available, and 
are served in the first-come-first-serve discipline. Agents join an agent queue upon acceptance of invitations if no customer is waiting, and serve customers in the order of their acceptance of invitations. 
Once a customer and an agent are matched, they leave the system immediately.   (Service times are irrelevant in our model.) 
Agent invitations are issued at the event times of customer arrivals and agents acceptance of invitations and at the independent Poisson event times driven by the difference of agent and customer queues. 
We assume that customer arrival process is Poisson and the invited agent response times are i.i.d. exponential.
(We discuss later how these assumptions can be relaxed.) 
The system state can be described by two variables. One tracks the difference of the agent queue and the customer queue. (Only one of those queues can be
positive at any time -- we assume that the non-idling condition is in force.)
Another variable is the target level of invited agents, which is also the actual number of invited agents.
(This is made possible by assuming that agent invitations can be issued or 
revoked 
instantaneously if needed.) Under our assumptions, the evolution
of the vector (`queue difference', `target level') is a CTMC.

Exact analysis of the model is prohibitively hard, and thus, we will analyze it in an asymptotic regime where the customer arrival rate becomes large while the distribution
of an agent response times is fixed. (Recall that service times do not matter here.)  This scaling regime is like the so-called many-server asymptotic regime.  In particular, the agent response process is modeled as an infinite-server model with an ``arrival" process controlled by the invitation scheme. 
Our main focus is the case when the customer arrival rate is constant. In this case, we study the process under the fluid and diffusion scalings. 
On the fluid scale, we show convergence to the fluid limit 
and uniform global stability of fluid limits (Theorem \ref{th-fluid-conv});
in addition we prove
the process stochastic stability and the limit-interchange
property (Theorem~\ref{th-fluid-dirac}) -- the sequence of fluid-scaled stationary
distributions converges to the distribution concentrated on a 
single point corresponding to zero queues.
The key technical challenge in the fluid-scale analysis stems 
from the fact that the target-level variable has to stay non-negative,
which creates a non-trivial boundary behavior.
Then, on the diffusion scale, we prove the convergence to the diffusion limit process (Theorem \ref{th-diff-conv})
and present the tightness and limit-interchange results (Theorem \ref{th-diff-tight}). (The latter can be obtained by adopting 
the approach in \cite{SY2012,St2013_tightness}. 
We give a high level sketch, but not all the lengthy details, in this paper,
because believe that the fluid-scale results are of greater interest and
importance for our model.)
When the customer arrival rate is time-varying, we give a fluid limit result (Theorem \ref{thm-fluidtv}).
 Our simulation experiments show good performance of the feedback scheme, as well as accuracy
 of the fluid limit approximation.

\subsection{Contributions and comparisons}

Our research contributes to workforce management in call centers, where both customers and agents must be managed properly in order to minimize operational costs \cite{AAM, GKM03}. There has been extensive study on the management of regular agents in the literature; however, there is a lack of stochastic models and analysis for the management of special agents/knowledge workers,  which can be engaged dynamically, in response to actual demand.  
A recently introduced (data-driven) Erlang ``S" model \cite{AFM} (more motivated from the management of regular agents), is related to ours in the 
general sense that it allows the dependence of agent-availability process on the customer queue length.  
A many-server queueing model with a random number of servers is studied in the context of virtual call centers or ride-sharing services in \cite{Ibrahim}. The number of agents available to work in a given period is modeled as an exogenous binomial random variable and the associated fluid model is studied to determine cost-minimizing staffing levels. In our model, the number of agents available is random but endogenous, depending upon the system state. 

Our model also relates to some extent to the literature on matching (double-ended) queues; see, e.g., \cite{GW, K66, LGK}.  In our system, however,
customer demand is matched by agents invited through a feedback control mechanism,
driven by the system state. This is very different from the standard matching (double-ended) queues since they assume that the entities to be matched arrive exogenously. 

It may appear that our model has some similarity with the classical make-to-stock (MTS)  queueing model in inventory theory; see, e.g., \cite{G94}. In the MTS model, demand of goods arrives as a Poisson process, and items are produced at a single-server or multi-server factory. 
A standard control algorithm for this system is as follows.
Upon the arrival of an order, a ``signal" is sent to the factory floor (FF) to produce one item, and an item is delivered to the customer from the finished goods (FG) inventory if it is available, and otherwise, the backlog is increased by one. This algorithm simply means that total inventory
(FF inventory plus FG inventory minus backlog) is kept at a constant level.
In comparison with our model, the backlog is like our ``customer" queue, the FG inventory is like our ``agent" queue, 
the net inventory level (FG inventory minus backlog) is like our 
``queue difference", and the FF inventory is like our ``pending agents".  
Our model and the algorithm differ from those for the MTS system in several aspects.\\ 
1) The key model difference is that we have infinite potential agent pool (``production capacity''), which is used to match any demand.\\
2) The ``invitation" schemes are  completely different, far beyond the fact that they apply to different models.
The fundamental difference is that in the MTS algorithm the total inventory is kept at a constant level,
which is computed as a function of the demand rate; our invitation scheme automatically adjusts to {\em any} demand rate (as long as it is 
sufficiently large) -- it does not need to be known or explicitly estimated. An ``analog"  of our scheme in the MTS context would be an algorithm
that automatically and dynamically adjusts to any demand rate.\\
3) Our results are, of course, different. In particular, they show that the invitation scheme is asymptotically optimal in the sense that
both customer and agent steady-state waiting times vanish as the system scale (customer arrival rate) becomes large. \\
4) We also note the robustness of our scheme. The asymptotic optimality is achieved for any setting of the algorithm parameters from a wide range, 
defined by certain simple conditions. In practice, it is easy to make sure that those conditions hold -- it suffices to have
only a rough idea about the system parameters (as opposed to exact or even approximate knowledge).\\
Finally, our model, almost as is,  might be useful for some MTS inventory systems, e.g., in settings where
there are many ``small" producers that may be activated and produce inventory after a delay. In such settings,
our asymptotic regime is natural.

As far as techniques are concerned, our work is relevant to the literature on the fluid and diffusion scale tightness and limit interchange in many-server asymptotic regime; see, e.g., \cite{GM08, GG13, GS12, SY2012,St2013_tightness} for an overview.  
 In our model, the process has a nontrivial boundary behavior. This presents new challenges for proving stability and fluid-scale limit-interchange. 
The two key difficulties are as follows.\\ 
 (a) Just because fluid limit trajectories have a unique stable point and satisfy a stable linear
ODE in the interior of the state space, 
does {\em not} directly imply uniqueness of trajectories and their uniform convergence 
to the stable point. The dynamics of trajectories on the boundary cannot be ignored.\\
(b) In the many-server regime in general, and our setting in particular,
the uniform convergence of fluid limit trajectories to the stable point
does {\em not} imply the system stochastic stability and/or fluid-scale tightness of stationary distributions.\\
To overcome these difficulties we develop a novel Lyapunov function and study the process on multiple time scales.
(See more detailed comments at the beginning of Section \ref{secProofFluid}.)
That is the main technical contribution of this paper.

\subsection{Organization of the paper}
The remainder of the paper is organized as follows. We will finish this section with basic notation and conventions below.  In section \ref{secModel}, we first describe the model and assumptions in detail. In section \ref{secMain},  we state the main fluid- and diffusion-scale results for the case
when the customer arrival rate is constant. 
 These results are proved in sections \ref{secProofFluid} and \ref{secProofDiff}, respectively.
Simulation experiments for the constant arrival rate case are provided in section \ref{secNumerical}. 
In section \ref{sec-timevarying}, we present a fluid limit result and simulations
for the time-varying arrival rate case. 
We will conclude and discuss future work in section \ref{sec-further-work}. 

\subsection{Basic notation and conventions}
\label{subsec-notation}

Sets of real and real non-negative numbers ($d$-dimensional vectors) are denoted by $\R$ and $\R_+$ ($\R^d$ and $\R_+^d$), respectively. $\NN$ denotes the set of positive integers.
The standard Euclidean norm of a vector $ x\in \R^n$ is denoted $\|x\|$. Vectors are viewed as row vectors.
For a vector $a$ or matrix $A$, we write their transposes as $a^T$ or $A^T$. 
We often write $x(\cdot)$ to mean the function (or random process) $(x(t),~t\ge 0)$.
 Let $D^k = D([0,\infty), \RR^k)$ 
  denote the space of $\RR^k$-valued functions defined on $[0,\infty)$ that are right continuous with left limits. 
For a real-valued function $x(\cdot): \R_+ \rightarrow \R$, we use either  $x'(t)$ or $(d/dt) x(t) $ to denote the derivative with respect to $t$, and for $x(\cdot): \R_+ \rightarrow \R^d$, we write $(d/dt) x(t) = (x'_1(t),..., x'_d(t))$. 
For a function $x(\cdot)$, we use $x(\infty)$ to denote its limit as $t \ra\infty$ if it exists. 
For a process $x(\cdot)$, we use $x(\infty)$ to denote a random variable having a steady state distribution of the process.
We use $\bone_A$ to denote an indicator function for a set $A$, where $\bone_A(x) = 1$ if $x\in A$ and  $\bone_A(x) = 0$ if $x\notin A$.
For a real number $x$, let $x^+ = \max\{x, 0\}$ and $x^{-} = - \min \{x, 0\}$, and let
$\sgn(x) = 1$ if $x >0$, $\sgn(x) =0$ if $x=0$ and $\sgn(x) =-1$ if $x <0$. The abbreviation  {\em w.p.$1$} means {\em  with probability $1$}.
Abbreviation {\em u.o.c.} means 
{\em uniform on compact sets} convergence of functions, with the argument (usually in 
$[0,\infty)$) determined by the context.
We write $x^r \to x \in \R^n$ to denote ordinary convergence in $\R^n$, and  $x^r \Rightarrow x$
to denote convergence in distribution of random variables taking values in space $\R^n$
equipped with the Borel $\sigma$-algebra.
Weak convergence of probability measures (on some Polish space) $\mu_n$ to $\mu$ will also be denoted as $\mu_n \Rightarrow \mu$. 
For a finite set of scalar functions $f_n(t), ~t\ge 0$, $n\in\NN$, a point $t$ is called
{\em regular} if for any subset $ \NN_o \subseteq \NN$ the 
derivatives
$$
\frac{d}{dt} \max_{n\in\NN_o} f_n(t) ~~\mbox{and}~~ 
\frac{d}{dt} \min_{n\in\NN_o} f_n(t)
$$
exist.   (To be precise, we require that each derivative is proper: both left and right derivatives exist and are equal.)
We use the familiar big-$O$ and small-$o$ notations for deterministic functions: for two real-valued functions $f$ and $g$, we write $f(x) = O(g(x))$ if $\limsup_{x\ra\infty} |f(x)/g(x)| < \infty$ and $f(x) = o(g(x))$ if $\limsup_{x\ra\infty} |f(x)/g(x)| =0$.

\section{Model and algorithm} \label{secModel}

\begin{figure}[htp]
\centering
\includegraphics[width=4.5in,height=2in]{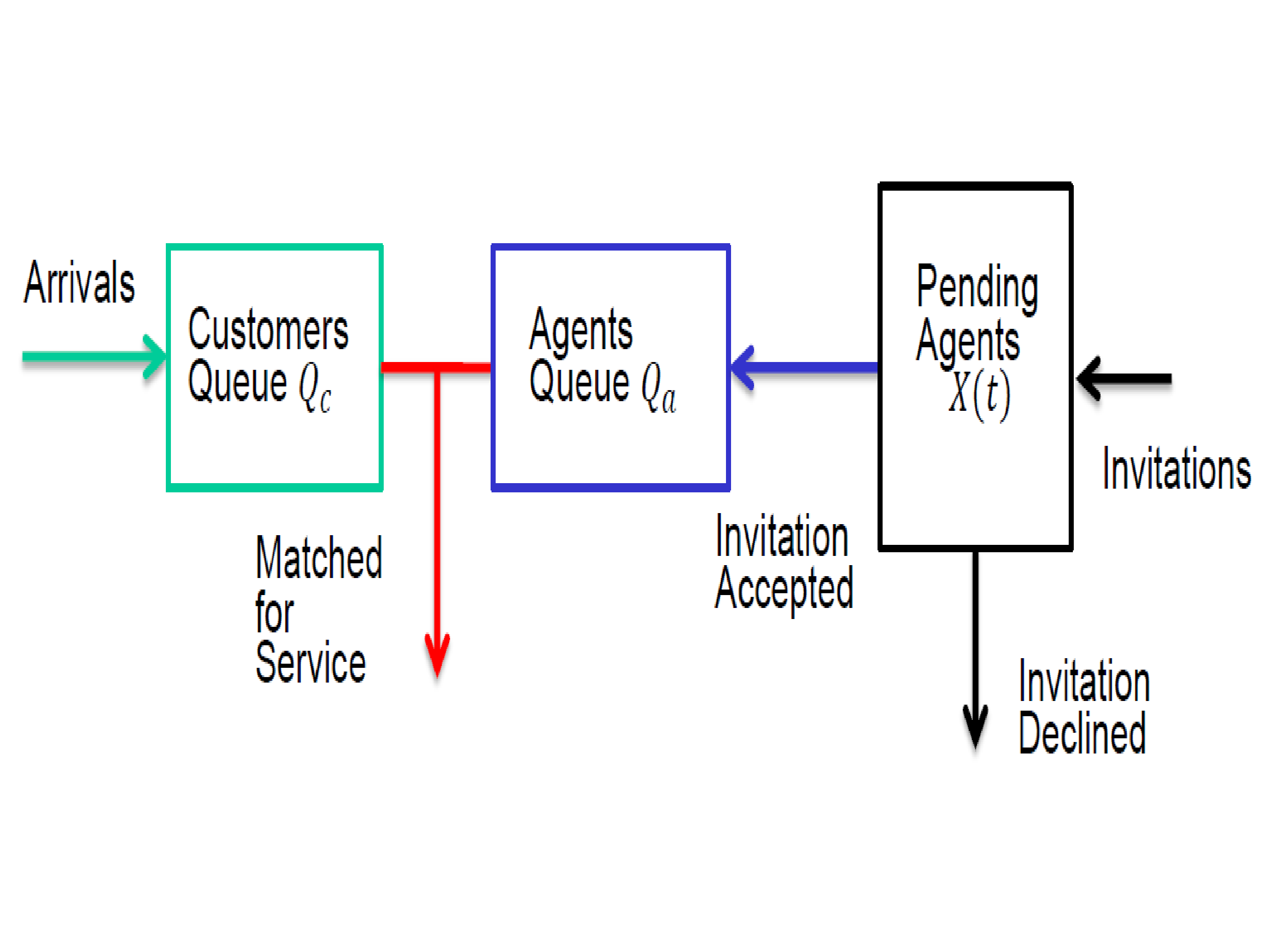}
\caption{An Agent Invitation System}
\label{fig1}
\end{figure}

Consider a customer contact center where agents are invited on demand and work at their distant ``homes". 
Customers arrive according to a Poisson process of rate $\Lambda >0$, and join a ``customer" queue waiting for an available agent and are served in the order of their arrival. 
Agents are invited to serve customers according to some scheme specified below. 
Once being invited, an agent will decide to accept or reject the invitation after some exponentially distributed random time. Let $\beta>0$ and $\tilde{\beta}>0$ be the rates at which an invited agent accepts or declines the invitation, respectively. 
Agents who accept their invitations will join an ``agent" queue waiting for a customer to arrive and serve customers in the order of their acceptance of the invitations. 
Once a customer and an agent are matched, they will leave the system simultaneously. This happens at the instant of either a customer arrival or an agent invitation acceptance.   Thus, we do not consider service times in this model.
Let $X(t)$ be the number of pending agents that have been invited but not decided to accept or decline the invitations at time $t$. 
 Let $Q_c(t)$ be the number of customers in the customer queue at time $t$ and $Q_a(t)$ be the number of agents in the agent queue at time $t$.   Define $Y(t) := Q_a(t) - Q_c(t)$, as the difference of the agent queue and customer queue at time $t$.
We assume that the non-idling condition holds, that is, agents do not idle when there are customers waiting in the customer queue, which implies that at each time $t$, either the customer queue or the agent queue must be empty. Figure \ref{fig1} depicts such an agent invitation system. 
 
The feedback invitation scheme in \cite{Patent2010}, let us label it as {\em Scheme A},
 works as follows. The scheme maintains a ``target" $X_{target}(t)$  for the number of invited agents $X(t)$.
The target $X_{target}(t)$ is changed by $\Delta X_{target}(t) = [-\gamma \Delta Y(t) - \ep Y(t) \Delta(t)]$ at each time $t$ when $Y(t)$ changes by $\Delta Y(t)$ (which can be either +1 or -1), where $\gamma>0$ and $\ep>0$ are the algorithm parameters and $\Delta(t)$ is  the time duration from the previous change of $Y$. 
 New agents are invited if and only if $X(t) < X_{target}(t)$, where $X(t)$ is the actual number of invited (pending) agents;
therefore, $X(t) \ge X_{target}(t)$ holds at all times. In addition, the target $X_{target}(t)$ is not allowed to go below zero, $X_{target}(t)\ge 0$;
i.e. if an update of $X_{target}(t)$ makes it negative, its value is immediately reset to $0$. Note that
$X_{target}(t)$ is {\em not} necessarily an integer.

To simplify our theoretical analysis, we consider a ``stylized" version of Scheme A, which has the same basic dynamics, but
keeps $X_{target}(t)$ integer and assumes that $X(t)=X_{target}(t)$ at all times; the latter is equivalent to assuming
that not only agent invitations can be issued instantly, but they can also be withdrawn at any time.
(In reality, it might not be feasible or desirable to withdraw the invitations. However, our simulations will confirm that
the performances of Scheme A and its stylized version are essentially same. See section \ref{secNumerical}.)
Given these assumptions, when pending agents decline invitations, it has no impact on the system state, 
because $X(t)$ is immediately ``replenished" by inviting another agent; therefore, in the analysis of stylized scheme, 
the events of declined invitations can be ignored (and parameter $\tilde \beta$ has no relevance).

Formally, the stylized scheme, which we label as {\em Scheme B}, is defined as follows. 
There are three types of mutually independent, and independent of the past, events that affect the dynamics of $X(t)$ and $Y(t)$
in a small time interval $[t,t+dt]$: a customer arrival with probability $\Lambda dt$, an agent acceptance with probability $\beta X(t) dt$, and
an additional event with probability $\ep |Y(t)| dt$. 

The changes at these event times are described as follows:\\
(i) Upon a customer arrival, $Y(t)$ changes by $\Delta Y(t) = -1$, and $X(t)$ increases by the average of $\gamma>0$.
For example, if $\gamma=2.6$ and $\Delta Y(t) = -1$,
then $\Delta X(t)=3$ with probability $0.6$ and $\Delta X(t) = 2$ with probability $0.4$.  To simplify the exposition, we assume 
from now on that $\gamma$ is an integer. 
\\
(ii)  Upon the acceptance of an invitation, it causes the change  $\Delta Y(t) = 1$ of $Y$, which in turn causes the change in $X(t)$ by $\Delta X(t)=-(\gamma \wedge X(t))$, that is, the change is by $-\gamma$ but $X(t)$ is kept to be nonnegative. 
\\
(iii) Upon the third type of event, if $X(t) \ge 1$, the change $\Delta X(t)=-\sgn(Y(t))$ occurs; and if $X(t)=0$, the change 
$\Delta X(t)=1$ occurs if $Y(t) < 0$ and $\Delta X(t)=0$ if $Y(t) \ge 0$.

All theoretical results in this paper concern Scheme B, which we consider in the rest of the paper, unless explicitly stated otherwise.
Under this scheme,
 the two-dimensional process $(Y, X)$ describing the system dynamics is a CTMC.  
We further assume that the parameters $\ep$, $\gamma$ and $\beta$ satisfy 
\beql{ep-assmp}
0 < \epsilon < \gamma^2 \beta / 4. 
\eeq

This is a technical condition. Its purpose is explained later, following equations \eqref{fluidxy3}--\eqref{fluidxyA}.

Informally, the scheme dynamics can be described as
$$
(d/dt) X = -\gamma (d/dt) Y - \epsilon Y, ~~ (d/dt) Y = \beta X - \Lambda,
$$
so that $X$ changes according to a negative feedback with respect to both the current value and increments of $Y$,
while $Y$ changes naturally, according to agents' invitation acceptances and new customer arrivals.

We remark that both Schemes A and B, of course, apply even when the customer arrival rate is time-varying. 
While our analysis is primarily focused  on the case of constant customer arrival rates (sections \ref{secMain}-\ref{secNumerical}),
we will provide some results for the time-varying case in section \ref{sec-timevarying}.
Also, the assumption of Poisson arrivals is not essential for our results and can be easily relaxed;
we believe the results can be generalized to allow non-exponential distribution of agent response times as well,
but this is less straightforward.


\section{Main Results} \label{secMain}

In this section, we state the main results of the paper. 
We consider a sequence of systems, indexed by the scaling parameter $r \in \RR_+$ and let $r\ra\infty$.  In the $r$-th system, the arrival rate is $\la r$,
while the parameters $\beta$,  $\ep$, $\gamma$ are constant.
The corresponding $r$-th process is $(Y^r, X^r)$,
where $Y^r=(Y^r(t), ~t\ge 0)$, $X^r=(X^r(t), ~t\ge 0)$.
 Define fluid-scaled processes with centering 
 \beql{xyf}
 (\bar{Y}^r, \bar{X}^r) := r^{-1} (Y^r, X^r - \la r /\beta) ,
 \eeq
 and diffusion-scaled processes 
 \beql{xyd}
  (\hat{Y}^r, \hat{X}^r): = \sqrt{r}  (\bar{Y}^r, \bar{X}^r) =  r^{-1/2} (Y^r, X^r - \la r /\beta) .
 \eeq

We first state a fluid limit result for $(\bar{Y}^r, \bar{X}^r)$ below. Note that for any $r$, $\bar{X}^r(t) \ge - \la /\beta$ for all $t$. 
\begin{thm} \label{th-fluid-conv} 
Consider a sequence of processes $(\bar{Y}^r, \bar{X}^r)$, $r\ra\infty$, with deterministic initial states such that
$ (\bar{Y}^r(0), \bar{X}^r(0)) \to (y(0),x(0))$ for some fixed $(y(0),x(0)) \in \RR^2$, $x(0) \ge -\la/\beta$.
Then,  these processes can be constructed on a common probability space, so that the following holds.
There exists a unique 
locally Lipschitz trajectory $(y,x)$,
such that, w.p.$1$, 
\beql{fluidconv}
 (\bar{Y}^r, \bar{X}^r) \to (y,x) \quad \text{u.o.c.}  \qasq r \ra\infty,
\eeq
where 
\beql{eq-fluid-bound}
x(t)\ge -\la/\beta, \quad t\ge 0, 
\end{equation}
and at any regular point $t\ge 0$ $($all points $t\ge 0$ are regular, except a subset of zero Lebesgue measure$)$,  the following holds: if $x(t) > -\la/\beta$,
\beqal{fluidxy}
y'(t) &=&\beta x(t), \non\\
x' (t)&=&  -\gamma \beta x(t) - \epsilon y(t) , 
\eeqa
and if $x(t) = -\la/\beta$,
\beqal{fluidxy-boundary}
y'(t) &=&-\lambda, \non\\
x' (t)&=&  [\gamma \la - \epsilon y(t) ] \vee 0. 
\eeqa
\end{thm}

The unique limit trajectory $(y,x)$ specified in Theorem~\ref{th-fluid-conv}
will be called a {\em fluid limit} starting from $(y(0),x(0))$.

Note that the second equation in \eqn{fluidxy} can also be written as  
\beql{fluidxy2}
x' (t)=  -\gamma y'(t) - \epsilon y(t).
\eeq
When the trajectory is away from the boundary,  the ODE \eqref{fluidxy} can be written as 
\beql{fluidxy3}
(d/dt) (y,x) = (y,x) A,
\eeq
where
\beql{fluidxyA}
A= \left[ \begin{array}{ll}
               0 & -\epsilon\\
                \beta & -\gamma \beta
                \end{array}
       \right] .
\eeq
The assumption  \eqref{ep-assmp} on the parameters guarantees that the matrix $A$ has two different negative eigenvalues. 
When the fluid limit trajectory hits the boundary $x(t)= - \la/\beta$ at some time $t$, $y$ will decrease at rate $\la$, that is, $y' = -\la$ until $y$ hits the value $\gamma \la /\ep$, and afterwards, the fluid limit will follow the trajectory of the ODE in \eqref{fluidxy} again. 
These observations imply, in particular, that 
the unique stable point of a fluid limit is $(0,0)$.

We then show the following stability and fluid-scale tightness results.

\begin{thm}
\label{th-fluid-dirac}
For all sufficiently large $r$, the system is stable, i.e.,  the Markov process $(Y^r, X^r)$ is positive recurrent. 
The sequence of 
stationary distributions of the fluid-scaled processes $(\bar{Y}^r, \bar{X}^r)$ converges to the Dirac measure concentrated at $(0,0)$.
\end{thm}

The following result describes transient behavior on the diffusion scale.

\begin{thm}
\label{th-diff-conv}
Suppose the sequence of deterministic initial states is such that
$(\hat{Y}^r(0), \hat{X}^r(0)) \to (\hat{Y}(0), \hat{X}(0))$,
where $(\hat{Y}(0), \hat{X}(0))$ is a fixed vector in $\RR^2$.
Then, 
\beql{fluidconv222}
 (\hat{Y}^r, \hat{X}^r) \RA (\hat{Y}, \hat{X})  \qasq r \ra\infty,
\eeq
where  $(\hat{Y}, \hat{X})$ is the unique solution to the SDE 
\beqal{diffxy}
\hat{Y}(t) &=& \hat{Y}(0) + \beta \int_0^t \hat{X} (s) ds  - \sqrt{2\la} W(t) , \non \\
\hat{X}(t) &=& \hat{X}(0) - \beta \gamma  \int_0^t \hat{X}(s) ds  - \ep \int_0^t \hat{Y} (s)ds  + \gamma \sqrt{2\la} W(t) , 
\eeqa
and $W$ is a standard Brownian motion.  
 The distribution of $(\hat{Y}(t), \hat{X}(t))$ is Gaussian with mean function $m(t)$
and covariance function $V(t)$ being unique solutions to the following ODEs, respectively:
\beqal{covxy}
\dot{m}(t) &=& m(t) A, \non \\
\dot{V}(t) &=& V(t) A + A^T V(t) + \sigma^T \sigma ,
\eeqa
where the matrix $A$ is given in \eqref{fluidxyA} and $\sigma := (- \sqrt{2\la}, \gamma \sqrt{2\la})$. The stationary distribution of $(\hat{Y}, \hat{X})$ is Gaussian with mean $(0,0)$ and covariance matrix
\beql{covxys}
V(\infty) = 
\left[ \begin{array}{ll}
               \frac{\la}{\beta \gamma} &  - \frac{\la}{\beta}\\
                - \frac{\la}{\beta} &  \frac{\la (\beta \gamma^2 + \ep)}{\beta^2 \gamma}
                \end{array}
       \right] .
\eeq
\end{thm}

Finally, we state the following diffusion-scale limit interchange result. 
We will not provide its detailed proof -- it can be done, using Theorems~\ref{th-fluid-dirac} and \ref{th-diff-conv} as the starting point,
and then following the approach given (for a different setting) in \cite{SY2012,St2013_tightness}.

\begin{thm}
\label{th-diff-tight}
The sequence
of stationary distributions of the diffusion-scaled processes $(\hat{Y}^r, \hat{X}^r)$  is tight. Consequently, 
given Theorem~\ref{th-diff-conv},
the limit-interchange holds:
the limit of stationary distributions of the diffusion-scaled processes $(\hat{Y}^r, \hat{X}^r)$ is equal to the stationary distribution of the limit diffusion process $(\hat{Y}, \hat{X})$.
\end{thm}

\section{Fluid scale analysis} 
\label{secProofFluid}

We first provide an overview of the fluid-scale analysis, which is the main part of this paper, and 
comment on the difficulties that arise.

The fluid limit $(y,x)$ in Theorem \ref{th-fluid-conv}
satisfies conditions \eqn{eq-fluid-bound}-\eqn{fluidxy-boundary}, including non-trivial boundary behavior
\eqn{fluidxy-boundary}. As a result, just because fluid limit trajectories satisfy a stable linear
ODE \eqn{fluidxy} (or, equivalently, \eqref{fluidxy3}-\eqref{fluidxyA}) away from the boundary, 
{\em does not directly imply uniqueness of trajectories and their uniform convergence 
to the stable point $(0,0)$.} Proving the latter properties is not straightforward; it requires studying the trajectory behavior on the boundary, in particular, how many times a trajectory can hit the boundary. 
To prove the uniform convergence of fluid limit trajectories
to the stable point, we develop a novel Lyapunov function using a norm associated with the eigenvectors of the matrix $A$ of the linear ODE  \eqref{fluidxy3}-\eqref{fluidxyA}.

The second major difficulty is that the convergence of fluid limit trajectories to the stable point $(0,0)$
{\em does not imply, in our setting, the system stochastic stability and/or fluid-scale tightness of stationary distributions}. To prove Theorem \ref{th-fluid-dirac} we need to employ a relatively involved analysis, 
considering first a sampled (discrete time) Markov chain and a Lyapunov-Foster criterion for its stability.
In turn, to show the negative drift of the Lyapunov function for the sampled chain,
we need to consider the original continuous time process (between sampled times)
on multiple time scales; this is because, roughly speaking, the scale of the  Lyapunov function decrease rate is very different on the state space boundary and away from it.

\subsection{Fluid models} \label{secPropFluid}


 We start studying properties of fluid limits by first considering 
{\em fluid models}, which are defined as locally Lipschitz continuous 
trajectories $(y,x)$, satisfying conditions \eqn{eq-fluid-bound}-\eqn{fluidxy-boundary}.
In other words, if a fluid limit exists (we do not even claim this yet), 
it is necessarily a fluid model. (The converse is
also true, as we will see, but we do not yet claim this either.)

\begin{lem}
\label{lem-fluid-model-prop}
For any initial state $(y(0),x(0))$, there is a unique
 fluid model starting from it. Moreover, uniformly
on the initial states from a given compact set,
$$
(y(t),x(t)) \to (0,0), ~~t\to \infty,
$$
and
$$
\max_{t\ge 0} \|(y(t),x(t))\| ~~\mbox{is bounded}.
$$
\end{lem}

Let us construct a fluid model starting from a given initial state $(y(0),x(0))$.
When the trajectory is away from the boundary,  the fluid model $(y,x)$ evolves
according to ODE \eqref{fluidxy3}-\eqref{fluidxyA}, 
and recall that the assumption on the parameters \eqref{ep-assmp} guarantees that there 
are two different negative eigenvalues of $A$, namely $-\infty < -\nu_2 < -\nu_1 <0$.

We choose the two corresponding eigenvectors to be $v_i=(\beta/\nu_i, -1)$, $i=1,2$.
Switching to the basis $v_1, v_2$, transforms any vector $u$ to $u B^{-1}$, where the matrix $B$ has 
$v_1$ and  $v_2$ as its first and second rows.
We will use norm $\|u\|_* = \|u B^{-1}\|$;
in other words, $\|u\|_* = (\alpha_1^2 + \alpha_2^2)^{1/2}$, where $\alpha_1$ and $\alpha_2$
are the coordinates of $u$ in the basis $v_1, v_2$, i.e.,  $u=\alpha_1 v_1 + \alpha_2 v_2$.

If we consider just the ODE  \eqref{fluidxy3} with initial state $(y(0),x(0)) =\alpha_1(0) v_1 + \alpha_2(0) v_2$, without any boundaries,
then the solution is
\beql{eq-ode-solution}
(y(t),x(t)) =\sum_{i=1,2} \alpha_i(0) e^{-\nu_i t} v_i .
\end{equation}
When/if the fluid model trajectory hits the boundary
$x=-\lambda/\beta$, say at time $t$, this can only happen if $y(t) \ge \gamma \lambda / \epsilon$ (recall the boundary dynamics in \eqref{fluidxy-boundary}); if so, then 
$x$ stays on the boundary ($x=-\lambda/\beta$, $x'=0$) and $y'=-\lambda$ until $y$ hits value $\gamma \lambda / \epsilon$;
at that time the fluid model starts obeying the ODE  \eqref{fluidxy3} again. Logically, it is possible that a fluid model trajectory
can hit the boundary multiple times. 
However, this cannot actually happen, as will follow from the following  Lemma~\ref{lem-norm-decrease}.

\begin{lem}
\label{lem-norm-decrease}
(i) There exists $\eta_1>0$ such that, for any fluid model, at any regular point $t$ such that $x(t) > -\lambda/\beta$,
\beql{eq-norm-decr1}
(d/dt) \|(y(t),x(t))\|_* \le -\eta_1 \|(y(t),x(t))\|_*.
\end{equation}
(ii) There exists $\eta_2>0$ such that, for any fluid model, at any regular point $t$ such that $x(t) = -\lambda/\beta$ and $y(t) \ge \gamma\lambda /\epsilon$,
\beql{eq-norm-decr2}
(d/dt) \|(y(t),x(t))\|_* \le -\eta_2.
\end{equation}
(iii) There exist $\eta>0$ and $C>0$ such that, for any fluid model, at any regular point $t$, $\|(y(t),x(t))\|_* \ge C$ implies
\beql{eq-norm-decr3}
(d/dt) \|(y(t),x(t))\|_* \le -\eta.
\end{equation}
\end{lem}


{\em Proof.} (i) This is the case when $(y(t),x(t))$ is away from  the boundary (and then satisfies the ODE  \eqref{fluidxy3}), 
and therefore satisfies the ODE  \eqref{fluidxy3}. Then,  the statement easily follows
from \eqn{eq-ode-solution}.

(ii) This is the case when a trajectory moves on the domain boundary, namely it has the form $(y(t),-\la/\beta) $ with $y(t)\ge \gamma\lambda /\epsilon$ and  $y'(t)=-\lambda$,
until it reaches point $(\gamma\lambda /\epsilon,-\lambda/\beta)$.
This also can be written as 
\beql{eq-new1}
(y(t),-\la/\beta) = \alpha(t) v_1 -(\alpha(t)-\la/\beta) v_2, 
\end{equation}
where
$$
\alpha(t) = \frac{y(t)-(\la/\beta) a_2}{a_1-a_2},
$$
and then 
\beql{eq-new2}
\alpha'(t)=-\lambda/(a_1-a_2). 
\end{equation}
We note that if $y(t)$ decreases from $+\infty$ to $\gamma\lambda /\epsilon$ (and $x(t)$ remains
equal to $-\la/\beta$), the corresponding value of $\alpha(t)$ must remain positive, and moreover, bounded away from $\lambda/\beta$,
i.e. $\alpha(t) \ge \lambda/\beta+\delta$ for $\delta>0$. (Otherwise, $(y(t),-\la/\beta)$ would become arbitrarily close to $(\la/\beta) v_1$,
which would imply $x'(t)>0$, as seen from \eqn{eq-ode-solution}. This is not possible, because the trajectory
cannot leave line $x= -\lambda/\beta$ before it reaches point $(\gamma\lambda /\epsilon,-\lambda/\beta)$.)
From this, \eqn{eq-new1} and \eqn{eq-new2}, we easily obtain \eqn{eq-norm-decr2}.

(iii) Time points $t$ where
$x(t) = -\la/\beta$ and $y(t) < \gamma\lambda /\epsilon$ -- this is the only condition not covered by (i) or (ii) -- cannot be regular.
Then, \eqn{eq-norm-decr3} is a corollary from (i) and (ii).
$\Box$

As a corollary of this lemma, we see that a solution to the ODE  \eqref{fluidxy3}, starting from point $(y(0),x(0))=(\gamma \lambda / \epsilon, - \lambda/\beta)$
cannot ever reach a point such that $x(t)= - \lambda/\beta$ and $y(t) \ge \gamma \lambda / \epsilon$.
(If not, the trajectory would then return to the point $(\gamma \lambda / \epsilon, - \lambda/\beta)$, which is impossible,
because according to  Lemma \ref{lem-norm-decrease} the norm $\|(y(t),x(t))\|_*$ is strictly decreasing, unless $(y(t),x(t)) \equiv (0,0))$.)
Therefore, the fluid model trajectory that we are constructing,
 will hit the boundary at most once. If that happens, it spends a finite time on the boundary,
reaches the point $(\gamma \lambda / \epsilon, - \lambda/\beta)$, and then follows the ODE \eqref{fluidxy3} thereafter. All claims of Lemma~\ref{lem-fluid-model-prop} easily follow.

\subsection{Proof of Theorem \ref{th-fluid-conv}} \label{secProofFluidconv}

Throughout this section, we are under the assumptions of Theorem \ref{th-fluid-conv}.
Given properties of the fluid models that we have already established,
in order to prove Theorem \ref{th-fluid-conv}, it suffices to show that w.p.1
from any subsequence of $r$ we can choose a further subsequence,
along which a u.o.c. convergence to a  fluid model holds.

Given the initial state $(Y^r(0), X^r(0))$, we construct
the processes $(Y^r, X^r)$, for all $r$, on the same probability space
via a common set of independent Poisson processes as follows: 
\beql{y1a}
Y^r(t) = Y^r(0) +N_2 \left(\beta \int_0^t  X^r(s) ds \right) - N_1(\la r t) ,
\eeq
\beql{x2a}
X^r(t) =Z^r(t) + \left(- \min_{0 \le s \le t} Z^r(s) \right) \vee 0, 
\eeq
\beqal{x1a}
Z^r(t) &=& X^r(0) +   \gamma N_1(\la r t)  - \gamma N_2 \left(\beta \int_0^t  X^r(s) ds \right)  \non\\
&& \qquad + N_3 \left(\ep \int_0^t \left(Y^r(s)\right)^{-} ds \right)  -  N_4 \left(\ep \int_0^t \left( Y^r(s) \right)^{+} ds \right),
\eeqa
and $N_i(\cdot)$, $i =1,...,4$, are mutually independent unit-rate Poisson processes. W.p.1, for any $r$, relations
\eqn{y1a}-\eqn{x1a} uniquely define the realization of $(Y^r, X^r)$ via the realizations
of the driving processes $N_i(\cdot)$. Relation \eqn{x2a} -- the ``reflection" at zero -- corresponds to the property that 
$X^r(t)$ cannot become negative.

The functional strong law of large numbers (FSLLN) holds for each Poisson process $N_i$:
\beql{eq-fslln}
N_i(rt)/r \to t, ~r\to\infty, ~~~\mbox{u.o.c., ~~w.p.1.}
\end{equation}

We consider the sequence of associated fluid-scaled processes $(\bar Y^r, \bar X^r)$  as defined in \eqref{xyf}. (Note that the processes $\bar X^r$ are centered.) Let a constant $m>\|(y(0),x(0))\|$ be fixed. For each $r$, on the same probability space as $(\bar Y^r, \bar X^r)$, let us define
a modified fluid-scaled process $(\bar{Y}_m^r, \bar{X}_m^r)$ as follows.
Let $(\bar{Y}_m^r, \bar{X}_m^r)$ start from the same initial state as $(\bar{Y}^r, \bar{X}^r)$, i.e., $(\bar{Y}_m^r(0), \bar{X}_m^r(0)) =(\bar{Y}^r(0), \bar{X}^r(0)) $. 
The modified process $(\bar{Y}_m^r, \bar{X}_m^r)$
follows the same path as $(\bar{Y}^r, \bar{X}^r)$ until the first time that $\|(\bar{Y}^r(t), \bar{X}^r(t))\| \ge m$.  Denote this time by  $\tau_m^r$.  
We then freeze the process $(\bar{Y}_m^r, \bar{X}_m^r)$  at the value $(\bar{Y}^r(\tau_m^r), \bar{X}^r(\tau_m^r))$, i.e.,  $(\bar{Y}_m^r(t), \bar{X}_m^r(t)) = (\bar{Y}^r(\tau_m^r), \bar{X}^r(\tau_m^r))$ for all $t \ge \tau_m^r$. 

The proof of the  convergence of fluid-scaled processes $(\bar{Y}^r, \bar{X}^r)$ will be in two steps, which are roughly as follows.
First, we show the convergence of $(\bar{Y}_m^r, \bar{X}_m^r)$ to a limit trajectory that 
behaves like a fluid model as long as the state norm is away from $m$.
(Here we will use the fact that the modified processes  $(\bar{Y}_m^r, \bar{X}_m^r)$ are uniformly bounded for all $r$ and $t\ge 0$ by construction.)
Second, for a given initial state, we choose the constant $m$ large enough,
so that the limit trajectory never reaches norm level $m$, and therefore 
it is the unique fluid model; this implies that
on any finite time interval, w.p.1, for all large $r$,
 $(\bar{Y}_m^r, \bar{X}_m^r)$  coincides with
$(\bar{Y}^r, \bar{X}^r)$, and therefore the latter converges to the fluid model.

\begin{lem} 
\label{lem-lip-m}
Fix $(y(0),x(0))$ and a finite constant $m>\|(y(0),x(0))\|$.
The following holds w.p.1.
From any subsequence of $r$,  we can find a further subsequence, along which
$(\bar{Y}_m^r, \bar{X}_m^r)$ converges u.o.c. to a Lipschitz continuous
trajectory $(y_m,x_m)$, which satisfies properties $\eqn{eq-fluid-bound}
- \eqn{fluidxy-boundary}$ at any regular time $t$ such that $\|(y_m(t),x_m(t))\| < m$.
\end{lem}

{\it Proof.} 
For the modified fluid-scaled processes  $(\bar{Y}^r_m, \bar{X}^r_m)$, we define the associated counting processes for upward and downward jumps: for $t \le \tau_m^r$, 
\begin{eqnarray} \label{mfp1}
\bar{Y}^{r, \uparrow}_m (t) &=& r^{-1} N_2 \left( r \beta \int_0^t [ \bar X^r_m(s)+\lambda/\beta] ds  \right) ,  \non\\
\bar{Y}^{r, \downarrow}_m (t) &=& r^{-1} N_1 (\la r t) , \non\\
\bar{X}^{r, \uparrow}_m (t) &=& r^{-1} \gamma N_1 (\la r t)  + r^{-1}  N_3\left( r \ep \int_0^t  \left(\bar Y^r_m(s)\right)^{-} d s \right) , \non\\
\bar{X}^{r, \downarrow}_m (t) &=& r^{-1} \gamma N_2 \left( r \beta \int_0^t [ \bar X^r_m(s)+\lambda/\beta] ds  \right) 
+ r^{-1}  N_4\left( r \ep \int_0^t  \left(\bar Y^r_m(s)\right)^{+} d s \right) ,
\end{eqnarray}
and for $t> \tau_m^r$, all these counting processes are frozen at their values at time $\tau_m^r$, that is, 
\beql{mfp1-frozen}
\bar{Y}^{r, \uparrow}_m (t) = \bar{Y}^{r, \uparrow}_m (\tau_m^r), \quad \bar{Y}^{r, \downarrow}_m (t)= \bar{Y}^{r, \downarrow}_m (\tau_m^r), \quad \bar{X}^{r, \uparrow}_m (t) = \bar{X}^{r, \uparrow}_m (\tau_m^r), \quad \bar{X}^{r, \downarrow}_m (t) = \bar{X}^{r, \downarrow}_m (\tau_m^r), \quad  t \ge \tau_m^r. 
\eeq

Using this notation, relations \eqn{y1a}-\eqn{x1a}, and the fact that for $0 \le t \le \tau^r_m $ the original and modified
processes, $(\bar{Y}^r, \bar{X}^r)$ and $(\bar{Y}^r_m, \bar{X}^r_m)$, coincide, we obtain for all $t\ge 0$:
\beql{mfp3}
\bar{Y}^r_m(t) = \bar{Y}^r(0) + \bar{Y}^{r, \uparrow}_m (t) - \bar{Y}^{r, \downarrow}_m (t),
\eeq
\beql{mfp4}
\bar{X}^r_m(t) = \bar{Z}^r_m(t) + \left( -\la/\beta - \min_{0\le s \le t}  \bar{Z}^r_m(s) \right) \vee 0 ,
\eeq
\beql{mfp41}
\bar{Z}^r_m(t) = \bar{X}^r(0) + \bar{X}^{r, \uparrow}_m (t) - \bar{X}^{r, \downarrow}_m (t).
\eeq

The counting processes $\bar{Y}^{r, \uparrow}_m, \bar{Y}^{r, \downarrow}_m, \bar{X}^{r, \uparrow}_m, \bar{X}^{r, \downarrow}_m$ are non-decreasing.
Using the FSLLN \eqn{eq-fslln} and the fact that the processes $\bar{Y}^{r}_m$ and $\bar{X}^{r}_m$ are uniformly bounded by construction,
we see that, w.p.1. for any subsequence of $r$, there exists a further subsequence along which the set 
of trajectories $(\bar{Y}^{r, \uparrow}_m,  \bar{Y}^{r, \downarrow}_m, \bar{X}^{r, \uparrow}_m, \bar{X}^{r, \downarrow}_m)$ converges u.o.c. to a 
set of 
non-decreasing Lipschitz continuous functions $(y^{\uparrow}_m, y^{\downarrow}_m, x^{\uparrow}_m, x^{\downarrow}_m)$.
But then the u.o.c. convergence of $(\bar Y_m^r, \bar X_m^r, \bar Z_m^r)$ to a set of Lipschitz continuous functions $(y_m, x_m, z_m)$ holds, where
\begin{eqnarray} \label{mfp5}
y_m(t) &=& y(0) + y^{\uparrow}_m(t) -  y^{\downarrow}_m(t), \\
x_m(t) &=& z_m(t) +  \left( -\la/\beta - \min_{0\le s \le t}  z_m (s) \right) \vee 0, 
\end{eqnarray}
\beql{mfp6}
z_m(t) = x(0) + x^{\uparrow}_m(t) -  x^{\downarrow}_m(t).
\eeq
Using this, and again the FSLLN \eqn{eq-fslln}, we can take the limit in \eqn{mfp1} to obtain:
\begin{eqnarray} \label{mfp7}
y_m^{\uparrow}(t) &=& \beta \int_0^t (x_m(s)+\la/\beta) ds, \non \\
y_m^{\downarrow}(t) &=&   \la t, \non \\
 x_m^{\uparrow} (t) &=&  \gamma \la t + \ep \int_0^t y_m^-(s) ds, \non \\
 x_m^{\downarrow}(t) &=&   \gamma \beta \int_0^t  (x_m(s)+ \la/\beta) ds + \ep \int_0^t y_m^+(s) ds,
\end{eqnarray}
It is easy to verify that properties $\eqn{eq-fluid-bound}
- \eqn{fluidxy-boundary}$ hold for the trajectory $(y_m, x_m)$.
 This completes the proof.
$\Box$



%


{\it Conclusion of the proof of Theorem \ref{th-fluid-conv}.}
For the given $(y(0),x(0))$, consider the corresponding (unique) fluid
model $(y,x)$. Let us choose $m > \max_{t\ge 0} \|(y(t),x(t))\|$.
Now let us apply Lemma~\ref{lem-lip-m}. W.p.1, from any subsequence of $r$ 
we can choose a further subsequence along which 
$(\bar{Y}_m^r, \bar{X}_m^r)$ converges u.o.c. to a Lipschitz continuous
trajectory $(y_m, x_m)$, which satisfies properties \eqn{eq-fluid-bound}
- \eqn{fluidxy-boundary}, as long as $\|(y_m(t),x_m(t))\| < m$.
But, as long as $\|(y_m(t),x_m(t))\| < m$, $(y_m, x_m)$ coincides with the fluid model
$(y,x)$. By the choice of $m$, this means that $(y_m, x_m)= (y,x)$.
Moreover, along the chosen subsequence, for any fixed $T>0$,
for all sufficiently large $r$, $(\bar{Y}_m^r, \bar{X}_m^r)$ and $(\bar{Y}^r, \bar{X}^r)$
coincide in the interval $[0,T]$. We see that, along the  chosen subsequence, u.o.c. convergence
of $(\bar{Y}^r, \bar{X}^r)$ 
to the fluid model $(y,x)$ holds. This means that w.p.1 the u.o.c. convergence
of $(\bar{Y}^r, \bar{X}^r)$ 
to $(y,x)$ holds for the original sequence of $r$. 
$\Box$

\subsection{Proof of Theorem \ref{th-fluid-dirac}} \label{secProofTight}

Given the uniform convergence of fluid limits in Lemma~\ref{lem-fluid-model-prop},
 to prove Theorem~\ref{th-fluid-dirac},  it suffices to prove the following

\begin{lem}
\label{lem-fluid-tight}
For all sufficiently large $r$, the process $(X^r,Y^r)$ $($and then $(\bar X^r, \bar Y^r)$$)$ is stable, with a unique stationary distribution.
The sequence of stationary distributions of  $(\bar X^r, \bar Y^r)$  is tight.
\end{lem}

Indeed, suppose Lemma~\ref{lem-fluid-tight} holds.
Consider stationary versions of the processes $(\bar Y^r(\cdot), \bar X^r(\cdot))$.
Fix arbitrary $\delta>0$ and then a sufficiently large compact set $B$ such that, uniformly in all (sufficiently large) $r$,
$\pr\{(\bar Y^r(0), \bar X^r(0)) \in B\} \ge 1-\delta$. 
Lemma~\ref{lem-fluid-model-prop} implies that we can choose a sufficiently large $T>0$ such that,
uniformly on all sufficiently large $r$ and initial states $(\bar Y^r(0), \bar X^r(0)) \in B$, we have
$$
\pr\{\|(\bar Y^r(T), \bar X^r(T))\| \le \delta ~|~ (\bar Y^r(0), \bar X^r(0))\} \ge 1-\delta.
$$
Therefore, for all large $r$, $\pr\{\|(\bar Y^r(T), \bar X^r(T))\| \le \delta \} \ge (1-\delta)^2$.
Since this is true for arbitrary $\delta>0$, we obtain the weak convergence of stationary distributions of $(\bar Y^r, \bar X^r)$
to the Dirac measure concentrated at $(0,0)$.

The approach we will take to prove Lemma~\ref{lem-fluid-tight}
is to first consider an embedded discrete-time Markov chain (DTMC), which is the original continuous-time chain sampled at a sequence of random stopping times, and show stability and the first moment bound for this DTMC  using a Lyapunov function drift criterion.
Specifically, we use the norm $\| \cdot \|_*$ defined in section \ref{secPropFluid} as the Lyapunov function.
We then use the relation between stationary distributions of the embedded and the original Markov chains.

{\em Proof of Lemma~\ref{lem-fluid-tight}.} Here we denote the fluid-scaled processes $s^r(t) = (\bar{Y}^r(t), \bar{X}^r(t))$, and to simplify the notation, we will drop the index  $r$, so $s(t) = (\bar{Y}(t), \bar{X}(t))$ below is the random process, not the fluid limit. 

We consider the embedded Markov chain. Fix constants $\delta>0$ and $\tau_{max}>0$.
For the process starting from a given state $s = s(0)$, consider the random stopping time $\tau_{\delta}(s)$, which is 
the first time $t$ when $|~\|s(t)\|_* - \|s\|_*~| \ge \delta$; we then define the stopping time
$\tau(s)=\tau_{\delta}(s) \wedge \tau_{max}$. 
Define a sequence of stopping times $\tau^{(k)}$, $k=1,2,...$ by
$$
\tau^{(1)} = \tau(s(0)), 
$$
$$
\tau^{(k)} = \tau^{(k-1)} +  \theta_{\tau^{(k-1)}} \tau^{(1)}, \quad k =2, 3, ...,
$$
where $\theta_{\cdot}$ is the random time shift operator associated with the process. In more detail,
$$
\tau^{(k)} = [\tau^{(k-1)} + \tau_{max}] \wedge \inf\{t> \tau^{(k-1)}: |~\|s(t)\|_* - \|s(\tau^{(k-1)}) \|_*~| \ge \delta\}, \quad k =2, 3, ...
$$

Consider the embedded discrete-time Markov chain $\hat s(k), ~k=0,1, \ldots$, using $\tau^{(k)}$ as sampling times. Specifically, if $s(t), ~t\ge 0$, is the original continuous time Markov process, then:
$$
\hat{s}(0) = s(0), \quad \hat{s}(k) = s(\tau^{(k)}), \quad k = 1,2,...
$$

Let $\Phi(s) =  \|s\|_*$. 
For the embedded chain $\hat{s}$, we show that, for some $C_1,C_2>0$, uniformly in $r$,
\beql{eq-key-drift}
\E [\Phi^2(\hat s(1)) - \Phi^2(\hat s(0)) ~|~\hat s(0)] \le -C_1  \Phi(\hat s(0)) + C_2.
\end{equation}
Note that $|\Phi(\hat s(1)) - \Phi(\hat s(0))|$ is uniformly bounded by the definition of $\tau(s)$.
Then, to prove \eqn{eq-key-drift} it suffices to show the following: for some constant $\delta_7>0$, for any sequence $r\to\infty$ and 
corresponding $\hat s(0)=\hat s^r(0)$ such that $\|\hat s^r(0)\|_* \uparrow \infty$, we have
\beql{eq-key-drift2}
\pr\{\Phi(\hat s^r(1)) - \Phi(\hat s^r(0))   \le -\delta_7\} \to 1.
\end{equation}
It suffices to consider a sequence such that the convergence
$$
\frac{1}{\|\hat s^r(0)\|_*} \hat s^r(0) \to \tilde s
$$
holds, for some vector $\tilde s$ with $\|\tilde s\|_*=1$.

We will study the behavior of the continuous time process $s(t)$, with initial state $s(0) = \hat s(0)$, 
on the interval $[0,\tau(s(0))]$. 
Before we proceed, we introduce some convenient (although somewhat abusive) notation.  For any vector $s=(y,x)$ we denote $s' = (y',x')$ where $y'=\beta x$, $x'= -\epsilon y - \gamma \beta x$;
in other words, these are the derivatives of the components of a fluid  trajectory $s(t)$ when $s(t)=s$.
Similarly, let $\|s\|'_*$ denote $(d/dt)\|s(t)\|_*$ when $s(t)=s$.

Suppose first that $\tilde s=(\tilde y, \tilde x)$ is such that $\tilde x > 0$.
Then, (the sequence of processes can be constructed on a common probability space, such that) w.p.1, u.o.c.
\beql{eq-slowdown-process}
s(t/\|s(0)\|_*) - s(0) \to \tilde s' t ~\mbox{and}~ \|s(t/\|s(0)\|_*)\|_* - \|s(0)\|_* \to \|\tilde s\|'_* t.
\end{equation}
From here \eqn{eq-key-drift2} follows. Indeed, we see that $\tau(s(0)) = \tau_{\delta}(s(0)) \to 0$, and therefore
\eqn{eq-key-drift2} holds with $\delta_7=\delta$.

Suppose now that $\tilde x =0$. Then, necessarily, $|\tilde y| >0$ and $\tilde y' =0$. If $\tilde y < 0$ then $\tilde x' > 0$, and this case
is treated the same way as the $\tilde x > 0$ case. Therefore, it remains to consider the case when 
$\tilde y >0$ and, consequently, $\tilde x' < 0$.

Consider the sub-case when
$$
[x(0) - (-\lambda/\beta)]/ |\tilde x' | \to \infty;
$$
then, we easily check that \eqn{eq-slowdown-process} still holds. This is the scenario when the time $\tau_{hit}$ for the $x(t)$ to hit boundary
$-\lambda/\beta$ is such that $\tau_{hit} \to 0$ and $\tau_{hit} \|s(0)\|_* \to \infty$; therefore, $\|s(t)\|_*$, which decreases at the rate 
$\|s(0)\|_* ~ \|\tilde s\|'_*$,
decreases by $\delta$ before time $\tau_{hit}$, and \eqn{eq-key-drift2} again follows.

Finally, consider the sub-case when (along a subsequence of $r$) 
$$
[x(0) - (-\lambda/\beta)]/ |\tilde x' | \to c \in [0,\infty).
$$
In this sub-case, $\tau_{hit} \|s(0)\|_* \to c$. Then, we consider the process such that in the interval $[0,\tau_{hit} \|s(0)\|_*]$ it is the process
with time slowdown, as in \eqn{eq-slowdown-process}, but from time $\tau_{hit} \|s(0)\|_*$ to infinity, the process continues in actual time, without
slowdown. W.p.1. in the limit we obtain the trajectory which satisfies \eqn{eq-slowdown-process} in the interval $[0,c]$, and then in the
interval $[c,\infty)$ we have $x(t)=-\lambda/\beta$ and $y'(t)=-\lambda$. In both intervals, the limit trajectory is such
that the norm $\|s(t)\|_*$ is decreasing at least at some
positive rate (see  Lemma~\ref{lem-norm-decrease}(iii)). We are done with proving \eqn{eq-key-drift2} and \eqn{eq-key-drift}.

From \eqn{eq-key-drift},  using standard Lyapunov-Foster argument (see, e.g., Chapter 13 in \cite{MT}), 
 we conclude that the embedded chain is stable for each sufficiently large $r$, and therefore has stationary distribution which is easily seen 
to be unique. Moreover, the stationary distributions are such that, uniformly in (sufficiently large) $r$,
\beql{eq-mean-bound-imb}
\E \Phi(\hat s(\infty)) \le C_2/C_1.
\end{equation}
We also observe that, for any fixed $C_3>0$, uniformly on all $\|s\|_*\le C_3$ and all $r$, $\E \tau(s) \ge C_4 >0$.
Let us choose $C_3$ large enough, so that for the embedded chain in steady-state, $\pr\{\|\hat s(\infty)\|_* \le C_3\} \ge 1/2$.

Now we use the relation between stationary distributions of the original continuous-time process and the sampled chain:
$$
\pr\{s(\infty) \in A\} = \frac{\E \left[\E \left[ \int_0^{\tau(s(0))} I\{s(t) \in A\} dt ~\big|~    s(0)=     \hat s(\infty) \right] \right]           } 
                                      {\E [ \tau(\hat s(\infty))] }
$$
(This relation is quite standard. For a proof, in a somewhat more general context, see Lemma 10.1 
 in \cite{St2003}.)
Then we see that our original continuous-time process is stable for each sufficiently large $r$,
and the stationary distributions are such that, uniformly in (sufficiently large) $r$, we have
\beql{eq-mean-bound}
\E \|s(\infty)\|_*  \le \frac{\E [\|\hat s\|_* + 2\delta] \tau_{max}} {C_4/2} \le C_5 \E \|\hat s(\infty)\|_* + C_6  \le C_7.
\end{equation}
The uniform bound \eqn{eq-mean-bound} on the expected norm in steady-state implies 
the tightness of stationary distributions.
$\Box$

\section{Diffusion scale analysis} \label{secProofDiff}

\subsection{Proof of Theorem \ref{th-diff-conv}} \label{secDiffConv}

We will use the following strong approximation of unit-rate Poisson processes (see Lemma 3.1 in \cite{K78}). 
\begin{lem} \label{lem-sa-p}
A unit-rate Poisson process $(\Pi(t): t\ge 0)$ can be realized on the same probability space as  a standard Brownian motion $(B(t): t\ge 0)$ such that the positive random variable $\xi$, given by
$$
\xi := \sup_{t\ge 0} \frac{|\Pi(t) - t - B(t)|}{ \log(2 \vee t)} < \infty
$$
has a finite moment generating function in a neighborhood of the origin. 
\end{lem} 

We will also need the following Lemma~\ref{lem-map}. Its proof follows from standard arguments (see, e.g., \cite{PTW}) and thus is omitted. 
\begin{lem} 
\label{lem-map}
Consider the mapping $\Phi: D^2 \ra D^2$  that takes  $(\phi_1, \phi_2) \in D^2$  into $(\psi_1, \psi_2) \in D^2$ determined by the integral representation: for each $t\ge 0$,
\beqal{mapping}
\psi_1(t) &=& \phi_1(t) + \beta \int_0^t \psi_2(s) ds  , \\
\psi_2(t) &=& \phi_2(t) - \beta \gamma \int_0^t \psi_2(s)ds - \ep \int_0^t \psi_1(s)ds,\non
\eeqa
where $\beta, \gamma, \ep$ are constants.
 Equations \eqref{mapping} determine $(\psi_1, \psi_2)$ uniquely. The mapping $\Phi$ is continuous in the topology of u.o.c. convergence.
\end{lem}

{\it Proof of Theorem \ref{th-diff-conv}.}
Since  $(\hat{Y}^r(0), \hat{X}^r(0)) \to (\hat{Y}(0), \hat{X}(0))$, we know that the fluid scaled processes $(\bar{Y}^r, \bar{X}^r)$  are such that $(\bar{Y}^r(0), \bar{X}^r(0)) \to (0,0)$. Then by Theorem~\ref{th-fluid-conv}, w.p.1,  $(\bar{Y}^r(t), \bar{X}^r(t))$ converges u.o.c. to the fluid limit  $(y(t),x(t)) \equiv(0,0)$ as $r\ra\infty$.  
Therefore,  w.p.1, 
\beql{diffpt}
\bar{Y}^r(t) \to 0 ,  \quad  \bar{X}^r(t) \to 0 , \quad    \int_0^t |\bar{Y}^r(s)| ds \to 0, \qandq \int_0^t | \bar{X}^r(s) | ds \to 0, ~~~\mbox{u.o.c.}
\eeq
This obviously implies that, w.p.1, on any finite time interval  $[0,T]$, the boundary of $\hat{X}^r$ is not hit for sufficiently large $r$,
i.e.  $\hat{X}^r(t) > -\lambda \sqrt{r}/\beta$. 
Using  representation \eqref{y1a} -- \eqref{x1a}  of the processes $(Y^r, X^r)$ (with $Z^r=X^r$)
and Lemma \ref{lem-sa-p}, after some manipulation, we see that 
there exist independent standard Brownian motions $B_i$, $i=1,2$, corresponding to the driving unit-rate Poisson processes $N_i$, $i=1,2$, all constructed on the same probability space, such that in the time interval $ t \in [0,T]$, w.p.1, 
\beqal{diffp1}
\hat{Y}^r (t)  &=& \hat{Y}^r (0) +  \int_0^t \beta \hat{X}^r(s) ds +  B_2\left(  \la t \right) -    B_1 (\la  t) +\Delta_1^r(t), \label{diffp1y}\\
\hat{X}^r  (t) &=& \hat{X}^r (0) - \gamma \int_0^t \beta \hat{X}^r(s) ds - \ep \int_0^t  \hat{Y}^r(s) ds +   \gamma B_1 (\la  t)   -   \gamma B_2\left(  \la t \right)    +  \Delta_2^r(t), \label{diffp1x}
\eeqa
where 
$$
\sup_{[0,T]} |\Delta_1^r(t)| = o(1), ~~~~ \sup_{[0,T]} |\Delta_2^r(t)|  =o(1).
$$
Letting $r\ra\infty$ and applying Lemma~\ref{lem-map} 
we obtain the probability 1, u.o.c., convergence of  $(\hat{Y}^r, \hat{X}^r)$  to a limit diffusion process $(\hat{Y}, \hat{X})$ which satisfies \eqref{diffxy}. This implies the claimed convergence in distribution.

Finally, the transient and stationary distributions of the limit diffusion process $(\hat{Y}, \hat{X})$  follow directly from linear SDEs; see Chapter 5.6 in \cite{KS96}.  This completes the proof of Theorem \ref{th-diff-conv}. $\Box$ 

\subsection{Comments on the proof of Theorem \ref{th-diff-tight}}
\label{secDiffTight}

This proof -- of the tightness of stationary distributions on diffusion scale -- can be obtained by
adapting the approach developed in \cite{SY2012, St2013_tightness}.  
That approach uses three steps. The first step is to show the fluid-scale ($r$-scale) tightness of the stationary distributions of the process;
in our context, it means proving that the stationary distributions of $r^{-1}(Y^r, X^r - \la r/\beta)$ are tight and asymptotically concentrate at $(0,0)$. We have done this step in Theorem \ref{th-fluid-dirac}.
The second step establishes the $r^{1/2+\kappa}$-scale tightness, for any $\kappa \in (0,1/2)$; namely,  tightness
of stationary distributions of   $r^{-1/2-\kappa}(Y^r, X^r - \la r /\beta) $. The argument of this step uses fluid-scale tightness
as a starting point, and follows that in \cite{SY2012}. The final step, is to show the diffusion-scale ($r^{1/2}$-scale) tightness;
it uses $r^{1/2+\kappa}$-scale tightness as a starting point, and follows the argument analogous to that in \cite{St2013_tightness}.

\section{Numerical Examples} \label{secNumerical}


In this section, we present some numerical examples to illustrate the behavior and performance 
of feedback agent invitation schemes, as well as accuracy of the approximations given by our theoretical results. 

Let us provide a general guidance on the setting of algorithm parameters.
For the (asymptotic) optimality of Scheme B, this setting does {\em not} need to be precise.
The rule of thumb is that $\epsilon$ should be ``as large as possible subject to condition $\epsilon < \beta \gamma^2/4$."
(This is because the smaller the $\epsilon$, the larger the system convergence time to the stationary point.)
For example, we can pick some reasonable value of $\gamma$, say $1$ or $2$.
Then, we can use some ``ballpark'' estimate for $\beta$. If this estimate is not very reliable, 
we set $\epsilon$ with large safety margin, say $\epsilon = (1/5)\beta \gamma^2/4$.
If the estimate of $\beta$ is considered more reliable, the choice of $\epsilon$ can be more aggressive,
say $\epsilon = (1/2)\beta \gamma^2/4$.

\begin{figure}[h]
\subfigure[Initial State $(Y(0), X(0)) = (0,0)$]{\label{plotfluid1}
\begin{minipage}{0.5\textwidth}
\centering
  \includegraphics[width=1\textwidth]{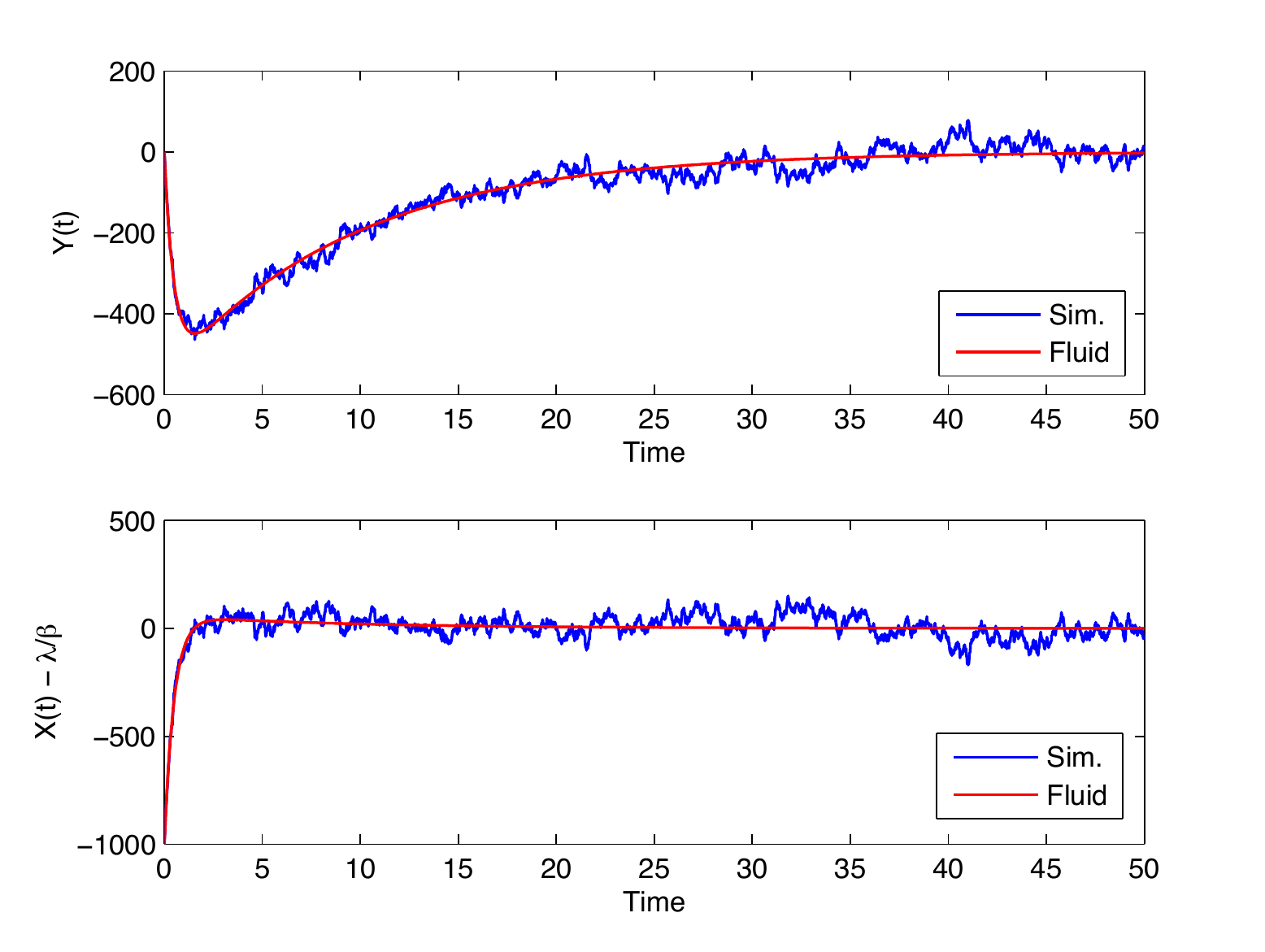}
\end{minipage}
}
\subfigure[Initial State $(Y(0), X(0)) = (1000,0)$]{\label{plotfluid2}
\begin{minipage}{0.5\textwidth}
\centering
  \includegraphics[width=1\textwidth]{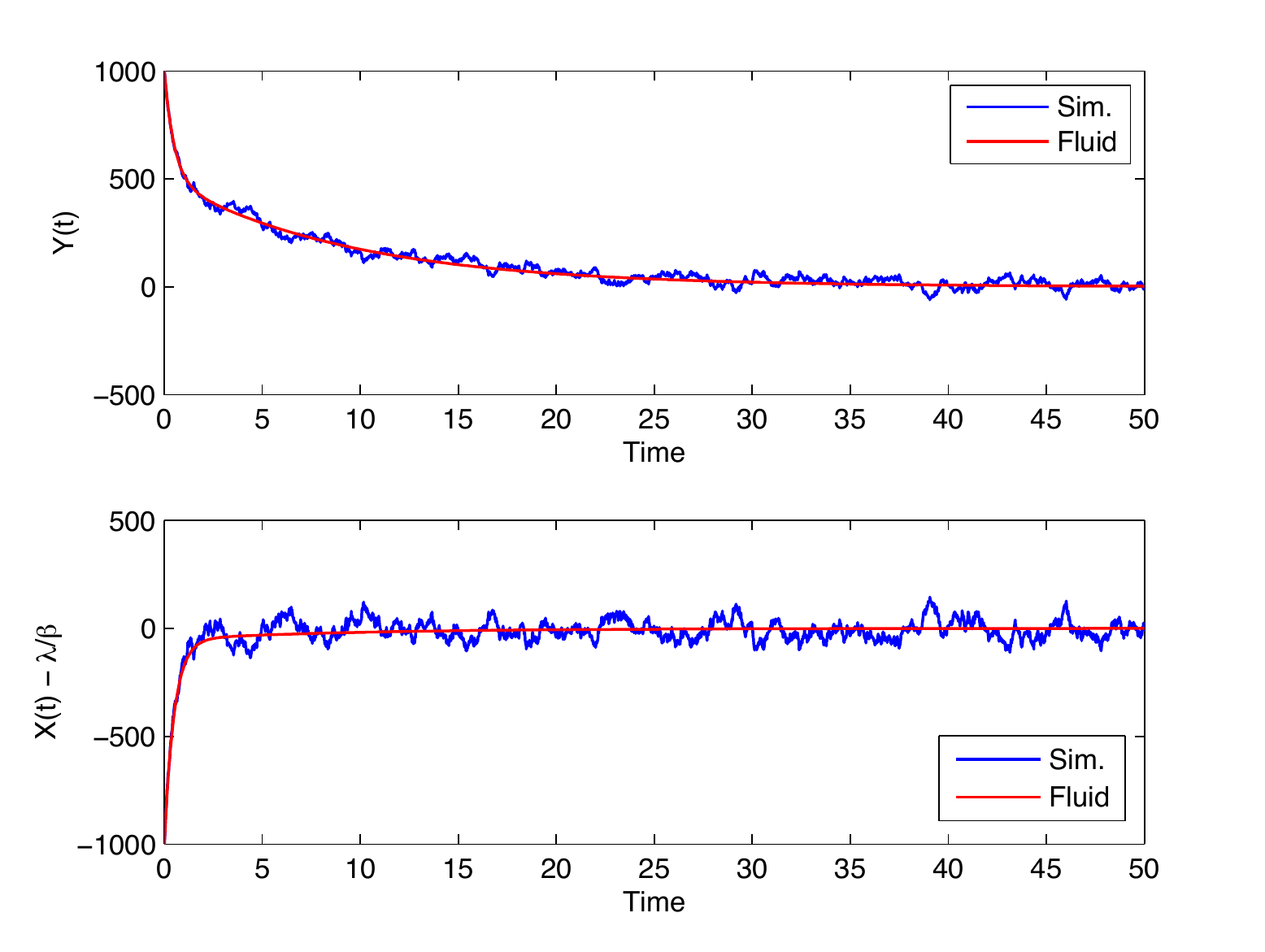}
\end{minipage}
}
\subfigure[Initial State $(Y(0), X(0)) = (0,2000)$]{\label{plotfluid3}
\begin{minipage}{0.5\textwidth}
\centering
  \includegraphics[width=1\textwidth]{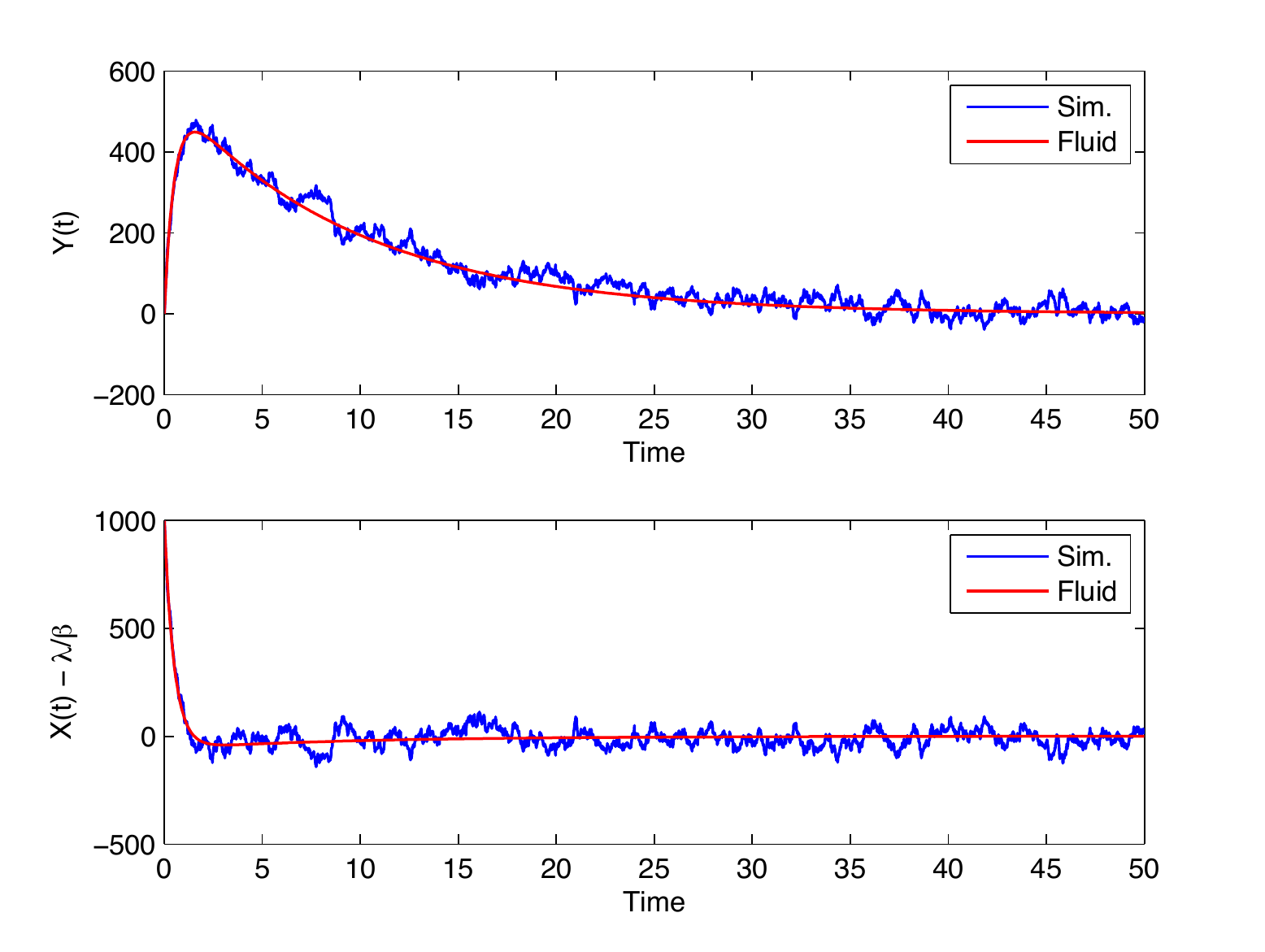}
\end{minipage}
}
\subfigure[Initial State $(Y(0), X(0)) = (-1000, 2000)$]{\label{plotfluid4}
\begin{minipage}{0.5\textwidth}
\centering
  \includegraphics[width=1\textwidth]{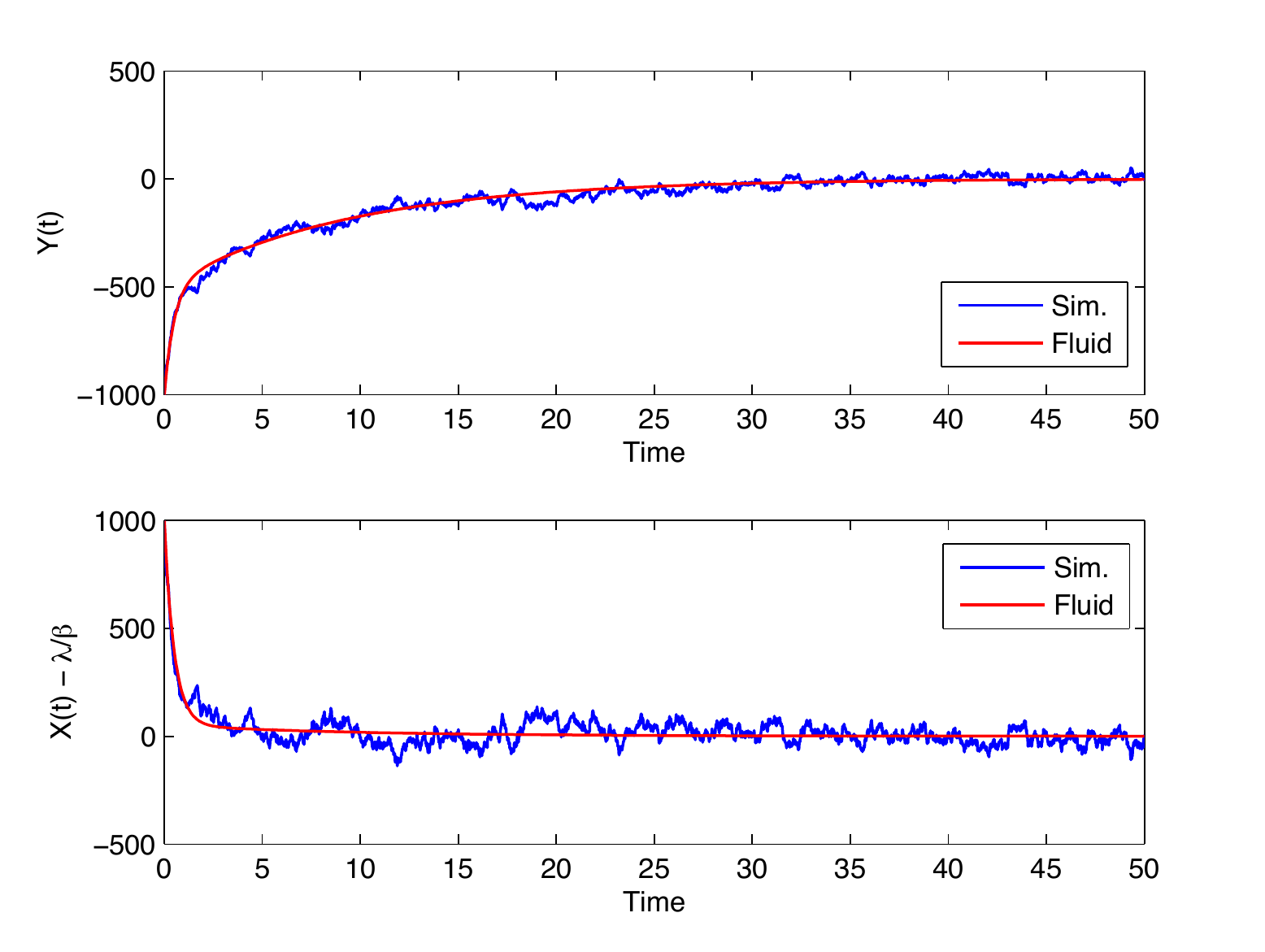}
\end{minipage}
}
\caption{Scheme B. Comparison of fluid approximations with simulations in one sample path. }
\label{plotfluid}
\end{figure}

First, we simulate Scheme B.
In the numerical examples, we use the following set of parameters: 
$$
\Lambda = 1000, \quad \beta = 1, 
\quad \gamma =2, \quad \ep = 0.2,
$$
We consider four initial conditions: 
(i) \ $(Y(0), X(0)) = (0,0)$;  (ii) \ $(Y(0), X(0)) = (1000, 0)$;  (iii) \ $(Y(0), X(0)) = (0, 2000)$; and
 (iv) \ $ (Y(0), X(0)) = (-1000, 2000)$.
In each case, we conduct  a simulation experiment of the system up to time 50 with these different initial conditions.
The comparisons in these cases are shown in Figure \ref{plotfluid}. 
We observe that the feedback scheme does bring the system (close) to its desired operating point, and fluid limit
provides a very good approximation of the system trajectory.

\begin{figure}[h]
\subfigure[ Fluid Model v.s. $X(t)$ and $X_{target}(t)$ ]{\label{plotfluidorig1}
\begin{minipage}{0.5\textwidth}
\centering
  \includegraphics[width=1\textwidth]{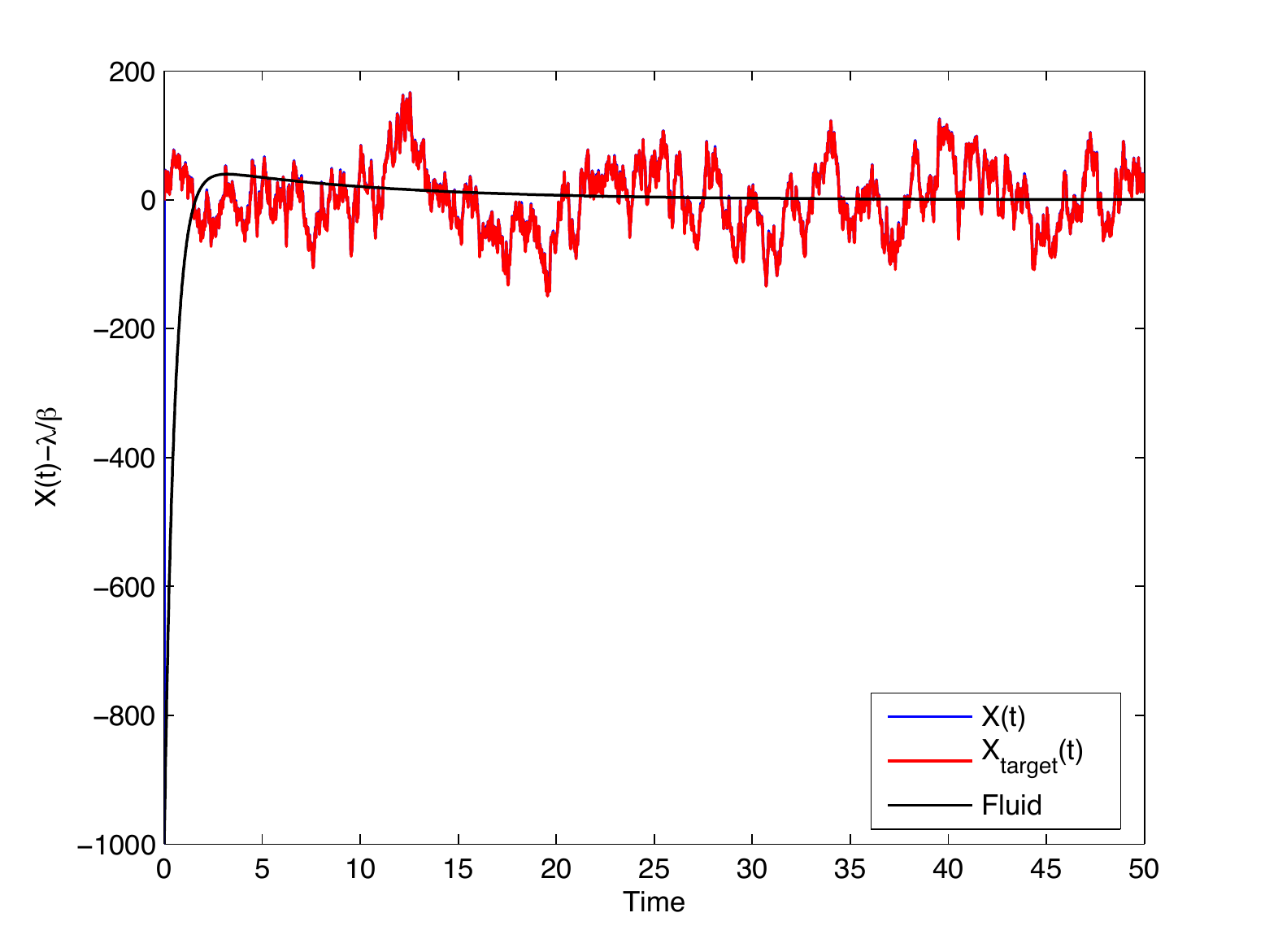}
\end{minipage}
}
\subfigure[ $X(t) - X_{target}(t)) $]{\label{plotfluidorig2}
\begin{minipage}{0.5\textwidth}
\centering
  \includegraphics[width=1\textwidth]{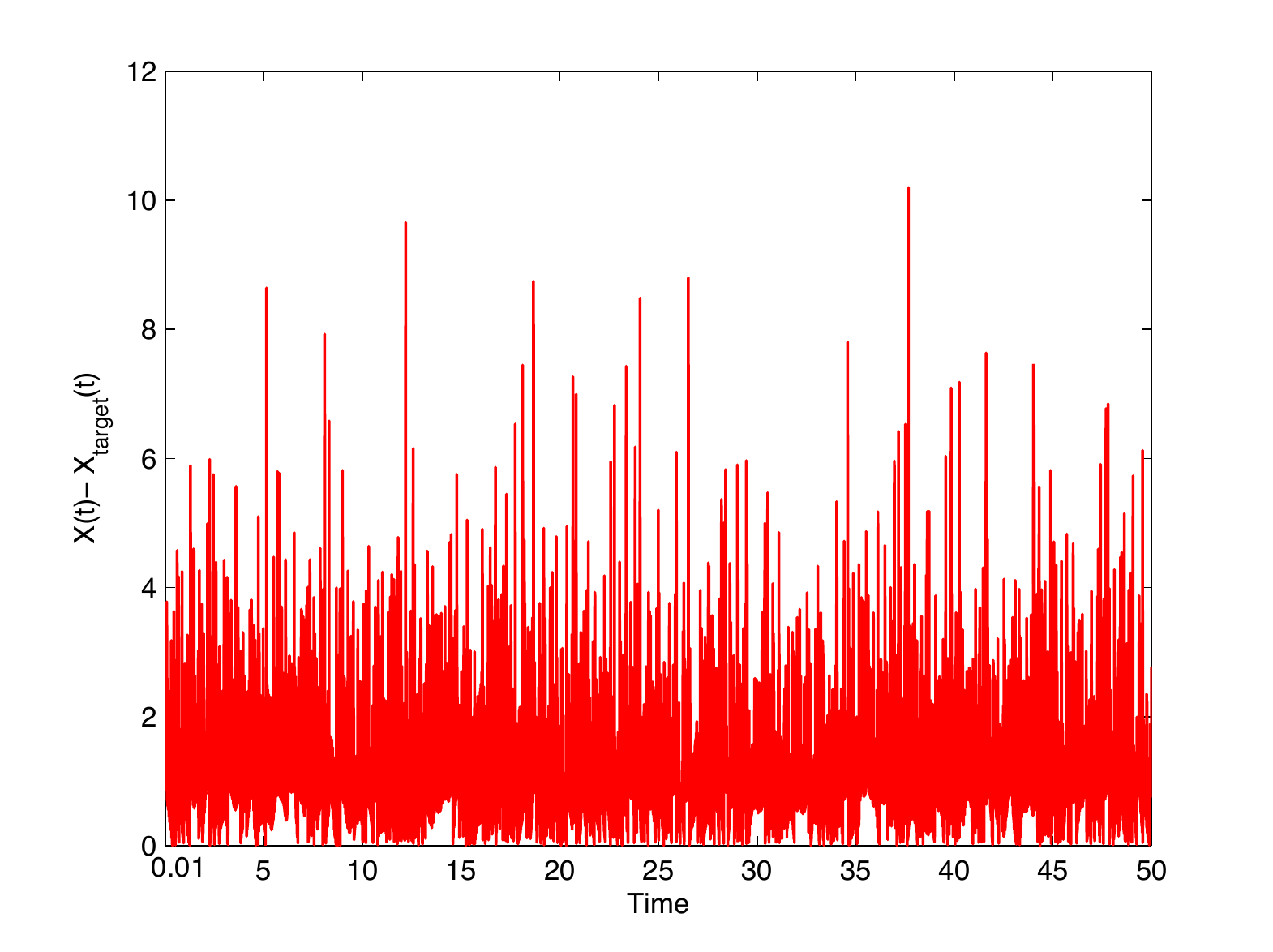}
\end{minipage}
}
\caption{Scheme A. }
\label{plotfluidorig}
\end{figure}

Second, 
we conduct a simulation experiment for Scheme A with the following parameter values 
$$
\Lambda = 1000, \quad \beta = 1, \quad \tilde{\beta} =1, \quad \gamma =2, \quad \ep = 0.2,
$$
 where $\tilde{\beta}$ is the rejection rate of invitations,
 and the initial state is $(Y(0), X(0), X_{target}(0)) =(0,0, 1000)$. 
The results are shown in Figure \ref{plotfluidorig1}. We see that the magnitude of the difference between $X_{target}$ and the actual number
of invited agents $X$ is very small (except at time 0) and can be regarded negligible compared to their scale.
This explains why the trajectories (of both $X_{target}$ and $X$) are well approximated
by the fluid trajectory, obtained for the Scheme B. This in turn provides a validation of Scheme B as an approximation of simpler and more practical
Scheme A.

\section{Time-Varying Customer Arrivals}
\label{sec-timevarying}

So far, we have considered systems in which customers arrive at a constant rate. 
In this section, we assume that customers arrive according to an inhomogeneous Poisson process. 
Both Scheme A and Scheme B introduced in section \ref{secModel} naturally work for time-varying customer arrivals.
Here we give a fluid limit result for Scheme B. 
We consider a sequence of systems $(Y^r, X^r)$, with the arrival rate function in the $r$-th system being
 $r \la(t), ~t\ge 0$, where $\la(\cdot)$ is a locally-bounded piece-wise continuous function with 
at most finite number of jumps on finite intervals.  Define fluid-scaled processes 
\beql{xyfv}
 (\tilde{Y}^r, \tilde{X}^r) := r^{-1} (Y^r, X^r).
\eeq
Note that the fluid-scaled process $\tilde{X}^r$ is not centered. 
\begin{thm} \label{thm-fluidtv}
Suppose that $ (\tilde{Y}^r(0), \tilde{X}^r(0)) \to (\tilde{y}(0),\tilde{x}(0))$ for some fixed $(\tilde{y}(0),\tilde{x}(0)) \in \RR^2$, $\tilde{x}(0) \ge 0$. 
Then the processes can be constructed on a common probability space, so that the following holds.
W.p.1, any subsequence of $r$ has a further subsequence, along which
\beql{fluidtvconv}
 (\tilde{Y}^r, \tilde{X}^r) \to (\tilde{y}, \tilde{x}) \quad \text{u.o.c.}  \qasq r \ra\infty,
\eeq
 where $(\tilde{y}, \tilde{x})$ is a  locally Lipschitz trajectory, such that
 at any regular point $t\ge 0$  the following holds: if $\tilde{x}(t) > 0$,
\beqal{fluidtvxy}
\tilde{y}'(t) &=&\beta \tilde{x}(t) - \la(t), \non\\
\tilde{x}'(t)&=&  \gamma \la(t) -\gamma \beta \tilde{x}(t) - \epsilon \tilde{y}(t) , 
\eeqa
and if $\tilde{x}(t) = 0$,
\beqal{fluidtvxy-boundary}
\tilde{y}'(t) &=&-\lambda(t), \non\\
\tilde{x}'(t)&=&  [\gamma \la(t) - \epsilon \tilde{y}(t) ] \vee 0. 
\eeqa

\end{thm}

Note that  the processes $(Y^r, X^r)$ can be represented from the equations \eqref{y1a} - \eqref{x1a} with $N_1(\la r t)$ being replaced by $N_1 \left(r \int_0^t \la(s)ds \right)$. 
 Given that, 
the proof of this theorem is a straightforward generalization of
the corresponding argument in the proof of Theorem \ref{th-fluid-conv}. We omit details.
Also note that Theorem~\ref{thm-fluidtv} does not claim uniqueness of the fluid limit trajectory.
However, the uniqueness easily follows in many cases of interest. For example, 
when $\la(\cdot)$ is piecewise constant (then the same argument as in the proof
of  Theorem \ref{th-fluid-conv} applies to each of the ``pieces"), or when the solution to \eqn{fluidtvxy} never hits the $\tilde{x}=0$
boundary.

\subsection{Numerical Example}

\begin{figure}[h]
\subfigure[ Initial State $(Y(0), X(0)) =(0,0)$  ]{\label{plotfluidtv1}
\begin{minipage}{0.5\textwidth}
\centering
  \includegraphics[width=1\textwidth]{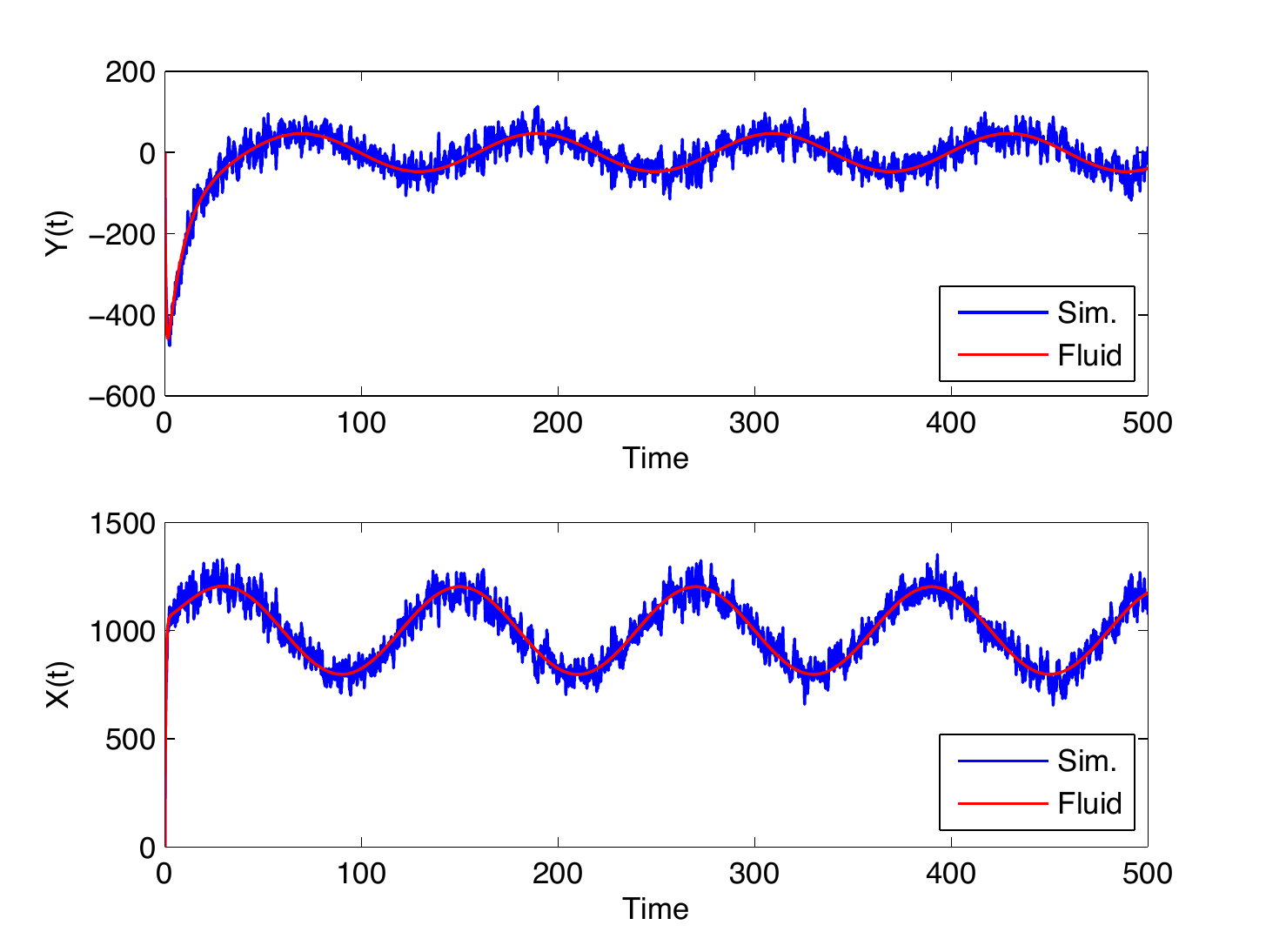}
\end{minipage}
}
\subfigure[  Initial State $(Y(0), X(0)) =(-1000, 2000)$]{\label{plotfluidtv2}
\begin{minipage}{0.5\textwidth}
\centering
  \includegraphics[width=1\textwidth]{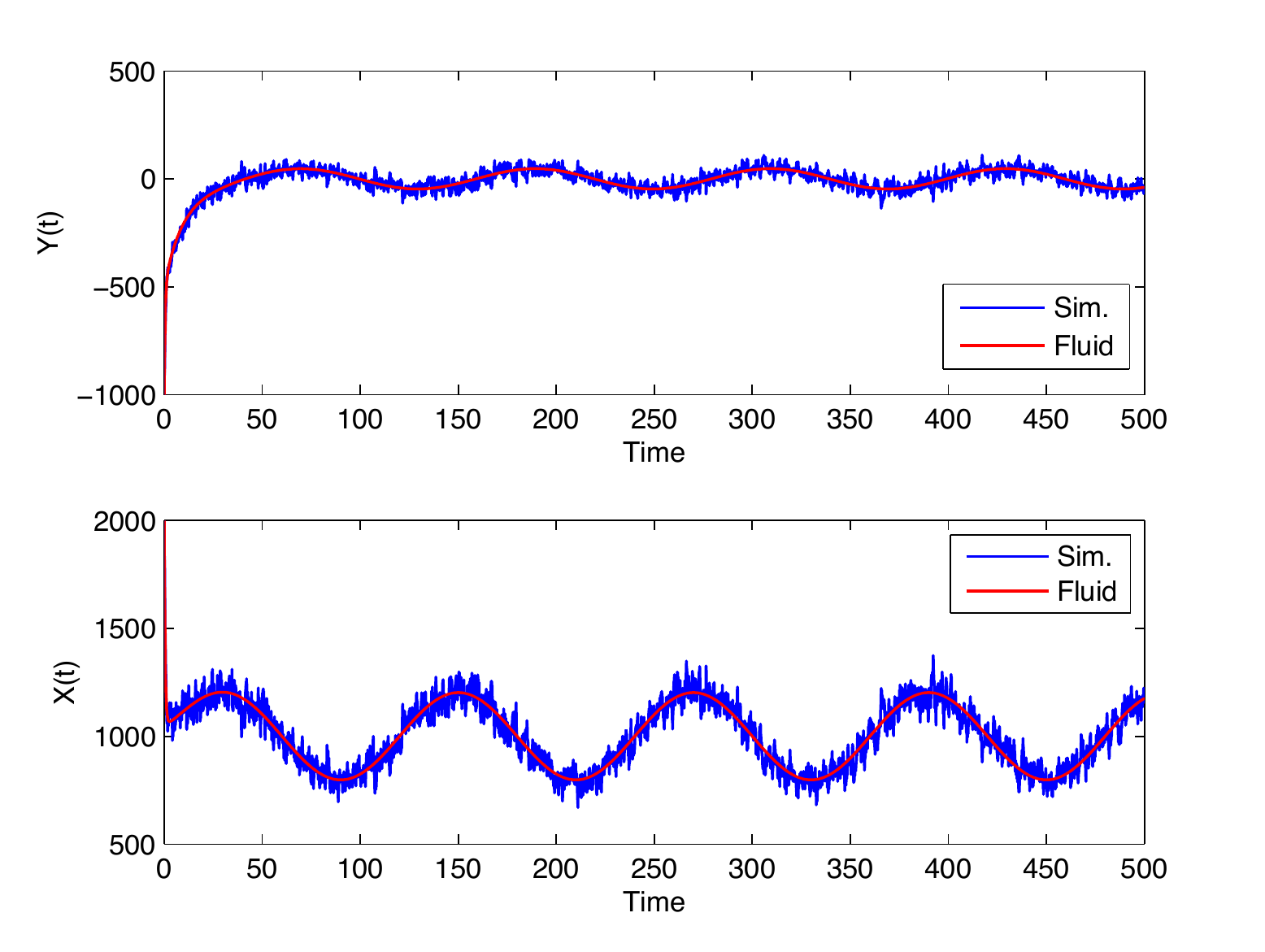}
\end{minipage}
}
\caption{Scheme B. Comparison of fluid approximations with simulations in one sample path when customer arrival rates are time-varying}
\label{plotfluidtv}
\end{figure}

We conduct a simulation experiment with the following set of parameters:
$$
\Lambda(t) = 1000+ 200 \sin(2\pi t/120), \quad \beta = 1, \quad \gamma =2, \quad \epsilon =0.2,
$$
and consider two initial conditions: (i) $(Y(0), X(0)) =(0,0)$ and (ii) $(Y(0), X(0)) =(-1000, 2000)$. In each case, we conduct a simulation experiment of the system up to time 500. See Figure \ref{plotfluidtv} for the comparisons. We observe that the fluid limit trajectory (which is unique in this case) provides a very good approximation of the system dynamics. We also observe that although the values of $Y$ do not converge to zero, they fluctuate around zero at a smaller scale than the system scale (comparing the magnitude of $Y$, less than 100, with that of the customer arrival rates, 1000). 

\section{Discussion and further work}
\label{sec-further-work}

{\bf Summary.}
We have proposed a new stochastic model
and the algorithm for a service system where agents are invited on demand.
The key feature of the algorithm is that it is very robust, requires essentially no knowledge of system parameters
and easily adapts to parameter changes. It is also very easy to implement.
To study its performance we prove the asymptotic results on the fluid and diffusion scales.
They, in particular, demonstrate the desirable performance of the scheme: both customer and agent waiting delays
vanish as the system scale increases to infinity.

{\bf Comment on the assumption of infinite number of available agents.}
In this paper we assume that the number of agents available to be invited is unlimited. 
It is a reasonable assumption, if the actual number of available agents, say $Mr$ in the system with parameter $r$, is much larger than 
the average number $(\lambda/\beta)r$ of agents required to meet the average demand rate $\lambda r$. 
Suppose that condition $M > \lambda/\beta$ does hold, but $M$ is finite.
With an infinite number of available agents we could and did ignore the fact that agents who are in the process of serving customers,
do not immediately return to the available agents' pool; with a finite number of agents, this can no longer be ignored, but suppose we keep this 
as a simplifying assumption (which approximately holds, for example, if service times are small). 
Then, our main results still hold, with essentially the same proofs. 
The only difference is that we will have an additional -- upper -- reflecting boundary  $ x \le M- \la/\beta$ for the for fluid-scaled number of invited agents;
the key  Lemma~\ref{lem-norm-decrease} easily generalizes to this case, and the entire analysis goes though as well.
In the case when we do {\em not} ignore the fact that agents performing service are unavailable to be invited, the system dynamics is more complicated.
Analyzing such a system may be a subject of future work.

{\bf Future work.}
Many more operational challenges for such systems remain open.
We have considered the model with single class of customers and single class of agents. 
 When there are multiple classes of customers and/or agent pools,
the design of an efficient invitation scheme presents new challenges. Finally,
in the model of this paper the service times are irrelevant; however, the service time cannot be ignored
the agent pools are finite. Exploring these and other related new problems may be a subject of future work.
In addition, our agent invitation scheme may be useful for the new generation of  {\em cloud-based} call centers. For instance, companies like Arise and LiveOps are providing platforms for businesses to implement their cloud-based call centers \cite{GLM, liveops, Arise2, Arise1}. The management of agents in a cloud-based system presents new operational challenges, in particular, how to guarantee the availability of qualified agents at any time. In some systems, the manager  sends invitation requests to qualified agents in order to meet the demand and satisfy the service level agreements. It will be interesting to investigate if our agent invitation scheme can be potentially used in cloud-based services in future work.  
Finally, in telemedicine operations, doctors and other medical specialists are valuable (and expensive) resources. Operational guidelines for telehealth services are provided  by the American Telemedicine Association \cite{ATA}.   The quality of care delivered via  telehealth systems requires timely interactions between health providers and patients - neither health providers nor patients would wait too long. 
Thus it will also be interesting to study in future work  if 
our ``agent" invitation scheme can be potentially useful for  telehealth management.

\end{document}